# BRIDGES AND NETWORKS: EXACT ASYMPTOTICS

By Robert D. Foley[1] and David R. McDonald[2]

*Georgia Institute of Technology and University of Ottawa*

We extend the Markov additive methodology developed in [*Ann. Appl. Probab.* **9** (1999) 110–145, *Ann. Appl. Probab.* **11** (2001) 596–607] to obtain the sharp asymptotics of the steady state probability of a queueing network when one of the nodes gets large. We focus on a new phenomenon we call a bridge. The bridge cases occur when the Markovian part of the *twisted* Markov additive process is one null recurrent or one transient, while the jitter cases treated in [*Ann. Appl. Probab.* **9** (1999) 110–145, *Ann. Appl. Probab.* **11** (2001) 596–607] occur when the Markovian part is (one) positive recurrent. The asymptotics of the steady state is an exponential times a polynomial term in the bridge case, but is purely exponential in the jitter case.

We apply this theory to a modified, stable, two node Jackson network where server two helps server one when server two is idle. We derive the sharp asymptotics of the steady state distribution of the number of customers queued at each node as the number of customers queued at the server one grows large. In so doing we get an intuitive understanding of the companion paper [*Ann. Appl. Probab.* **15** (2005) 519–541] which gives a large deviation analysis of this problem using the flat boundary theory in the book by Shwartz and Weiss. Unlike the (unscaled) large deviation path of a Jackson network which jitters along the boundary, the unscaled large deviation path of the modified network tries to avoid the boundary where server two helps server one (and forms a bridge). In the fluid limit this bridge does collapse to a straight line, but the proportion of time spent on the flat boundary tends to zero.

This bridge phenomenon is ubiquitous. We also treated the bathroom problem described in the Shwartz and Weiss book and found the bridge case is present. Here we derive the sharp asymptotics of the steady state of the bridge case and we obtain the results consistent

Received May 2002; revised June 2004.
[1]Supported in part by NSF Grant DMI-99-08161.
[2]Supported in part by NSERC Grant A4551.

with those obtained in [*SIAM J. Appl. Math.* (1984) **44** 1041–1053] using complex variable methods.

**1. Introduction.** In Section 2 we analyze the spectral radius of the Feynman–Kac transform of a Markov additive process. In Section 5 we develop a ratio limit theorem for Markov additive processes. In Section 6 we use these results to extend the method in [8, 14] to a general theory for the sharp asymptotics of the steady state probability $\pi$ of queueing networks when the queue at one node gets large. Finally, in Section 7 we apply these general results to obtain the asymptotics of $\pi(\ell, y)$, the steady state of a two node modified Jackson network, as $\ell \to \infty$.

A two node Jackson queueing network is described in the companion paper [9] which we refer to as Part I. The servers at node two, when idle, can assist the server at node one. Allowing one of the servers to help can completely change the behavior of the network. Clearly, a large deviation where the first node gets large tends to avoid emptying the second node because the idle node will help the overloaded node and make the large deviation less likely. The large deviation path to a point where the queue at node one is of size $\ell$ and the queue at node two is empty or of fixed size will tend to be a bridge. Asymptotically the large deviation path spends no time at states where the second queue is idle, as was seen in Part I.

This paper explains the dichotomy between the positive recurrent (jitter case) when the idle node does not help much and the large deviation path jitters along the axis where the second queue is idle and the transient or null recurrent (bridge case) when the idle node helps a lot giving rise to the bridge behavior. The former gives an exponential decay of $\pi(\ell, y)$ as $\ell \to \infty$, while the latter has an additional polynomial factor.

The bridge phenomenon is ubiquitous. In Section 8 we use this theory to revisit the *bathroom problem* discussed in [6, 14, 22]. Depending on the parameters, we obtain the jitter case discussed in [14] or a bridge case. In Section 8 we work out the asymptotics of the bridge case and we find exactly the same asymptotics of the steady state probability as was obtained in [6] using complex variable methods.

There are several approaches to obtaining the sharp asymptotics of the steady state probability of a queueing network. The oldest is the exact solution in product form when the network is quasi-reversible. The book by Serfozo [21] gives the state of the art. These product form solutions are a minor miracle, but they are very fragile as the slightest change in network dynamics can destroy the product form. The compensation method developed by Adan [1] is a generalization where one represents the steady state probability as an infinite linear combination of product measures. This method does allow one to attack nonproduct form steady states, but it is essentially a two-dimensional theory. It is hard to derive the asymptotics of the steady state and this is generally the most useful quantity. It is also hard to see how this theory could handle bridges where the steady state probability decays like a polynomial times an exponential.

The complex variable method used in [6] reduces the functional equation satisfied by the two-dimensional $z$-transform of the steady state probability to a Riemann Hilbert boundary value problem. One obtains a representation for the $z$-transform, which, in principle, determines the steady state, but again it is essentially a two-dimensional theory and it can be hard to derive useful quantities like the asymptotic behavior of the steady state probabilities. We were unable to derive the asymptotics of the solution of the coupled processors model (a special case of our modified Jackson network) obtained by complex variable methods in [4]. Fortunately for us, [6] does provide asymptotics and they agree with ours!

Finally, the matrix geometric method [16] has been recently extended in [23] to obtain the asymptotics of quasi birth and death QBD processes with infinite phase. In fact, these processes are Markov additive and closer inspection reveals a close parallel with the results in [8] and [14]. We discuss this parallel in Section 8.

In future work we will investigate the sharp asymptotics of the mean time until the queue at one node of the modified Jackson network overloads, as well as the Yaglom limit of the distribution of the queue at the second node.

**2. Markov additive processes and the Feynman–Kac transform.** In this section we define the Feynmac–Kac transform $J_\gamma$ of the transition kernel $J$ of a Markov additive process. The Markovian part of $J_\gamma$ will be denoted by $\hat{J}_\gamma$. Lemma 1 gives a representation of the spectral radius of $\hat{J}_\gamma$ and a condition for determining whether $\hat{J}_\gamma$ is $R$-recurrent or $R$-transient, where $1/R$ is the spectral radius of $\hat{J}_\gamma$.

Let $(V, Z) \equiv (V[n], Z[n]), n = 0, 1, 2, \ldots$, be a Markov additive process (see [17]) on a countable state space $\mathbb{Z} \times \hat{S}$, where $\mathbb{Z} = \{\ldots, -1, 0, 1, \ldots\}$. Given any state $x$, we denote the first component by $x_1$ and the second component by $\hat{x} \in \hat{S}$. The process $(V, Z)$ is a Markov chain with the following additional structure:

$$\begin{aligned} \Pr\{(V[n+1], Z[n+1]) = x | (V[n], Z[n]) = z\} \\ = J(z, x) \\ = J((z_1, \hat{z}), (x_1, \hat{x})) \\ = J((0, \hat{z}), (x_1 - z_1, \hat{x})), \end{aligned}$$

where $J$ denotes the transition kernel of $(V, Z)$. Note that the marginal process $Z$ forms a Markov chain with state space $\hat{S}$. We assume that $Z$ is irreducible, and we let $\hat{J}$ denote the transition kernel of $Z$.

As in [17], let $J_\gamma$ denote the *Feynman–Kac* transform of $J$. That is, for any real $\gamma$,

$$J_\gamma((x_1, \hat{x}), (y_1, \hat{y})) = e^{\gamma(y_1 - x_1)} J((x_1, \hat{x}), (y_1, \hat{y})). \quad (1)$$

Let $\hat{J}_\gamma$ be the Markovian part of $J_\gamma$, that is,

$$\hat{J}_\gamma(\hat{x}, \hat{y}) = \sum_v J_\gamma((0, \hat{x}), (v, \hat{y})). \quad (2)$$

Since $\hat{J}$ is irreducible, we know that $\hat{J}_\gamma$ is irreducible. We assume $\hat{J}_\gamma(\hat{x}, \hat{y}) < \infty$ for all $\hat{x}, \hat{y} \in \hat{S}$. Let $r(\hat{J}_\gamma)$ be the spectral radius of $\hat{J}_\gamma$ and let $R(\hat{J}_\gamma) = 1/r(\hat{J}_\gamma)$ be the associated convergence parameter. By Lemma 2 in [10], $\Lambda(\gamma) = \log(r(\hat{J}_\gamma))$ is a closed convex proper function [and so is $r(\hat{J}_\gamma)$]. As an aside, the next paragraph gives a short proof that $\Lambda(\gamma)$ is a convex function of $\gamma$.

First, for every $\hat{y}$, $A_n(\gamma) := (\hat{J}^n_\gamma(\hat{y}, \hat{y}))^{1/n} \to r(\hat{J}_\gamma)$. Next, for any $n$, we verify the function $\log(A_n(\gamma))$ is convex in $\gamma$ by checking that the second derivative is nonnegative. The second derivative is $A''_n/A_n - (A'_n)^2/A_n^2$. By calculation,

$$\frac{(A'_n)^2}{A_n^2} = \frac{1}{n^2} \left( \sum_v v \frac{1}{\hat{J}^n_\gamma(\hat{y}, \hat{y})} J^n_\gamma((0, \hat{y}), (v, \hat{y})) \right)^2$$

and

$$\frac{A''_n}{A_n} = \frac{1}{n} \left( \sum_v v^2 \frac{1}{\hat{J}^n_\gamma(\hat{y}, \hat{y})} J^n_\gamma((0, \hat{y}), (v, \hat{y})) \right.$$
$$\left. - \left(1 - \frac{1}{n}\right) \left( \sum_v v \frac{1}{\hat{J}^n_\gamma(\hat{y}, \hat{y})} J^n_\gamma((0, \hat{y}), (v, \hat{y})) \right)^2 \right).$$

Hence,

$$\frac{d^2}{d\gamma^2} \log(A_n(\gamma)) = \frac{1}{n} \left( \sum_v v^2 \frac{1}{\hat{J}^n_\gamma(\hat{y}, \hat{y})} J^n_\gamma((0, \hat{y}), (v, \hat{y})) \right.$$
$$\left. - \left( \sum_v v \frac{1}{\hat{J}^n_\gamma(\hat{y}, \hat{y})} J^n_\gamma((0, \hat{y}), (v, \hat{y})) \right)^2 \right).$$

The quantity in brackets above is equal to the variance of the measure $J^n_\gamma((0, \hat{y}), (\cdot, \hat{y}))/\hat{J}^n_\gamma(\hat{y}, \hat{y})$, which must be nonnegative. Since a limit of convex functions is convex, it follows that $\Lambda(\gamma)$ is a convex function of $\gamma$ [and so is $r(\hat{J}_\gamma)$].

Let

$$f_\gamma^n(\hat{y},\hat{y}) = \sum_{\hat{y}[1]\neq\hat{y},\ldots,\hat{y}[n-1]\neq\hat{y},\hat{y}[n]=\hat{y}} \hat{J}_\gamma(\hat{y},\hat{y}[1])\hat{J}_\gamma(\hat{y}[1],\hat{y}[2]) \times \cdots$$
$$\times \hat{J}_\gamma(\hat{y}[n-2],\hat{y}[n-1])\hat{J}_\gamma(\hat{y}[n-1],\hat{y}[n])$$

and define the transform $\Psi(\gamma,u) = \sum_n f_\gamma^n(\hat{y},\hat{y})u^n$. Note that $\Psi(\gamma,e^{-\zeta}) = \psi(\gamma,\zeta)$, defined at (3.2) in [17]. In [17], it is shown that if the set $\mathcal{W} = \{(\gamma,\zeta); \psi(\gamma,\zeta) < \infty\}$ is open and if there exists a $\Lambda(\gamma) < \infty$ such that $\psi(\gamma,\Lambda(\gamma)) = 1$, then $\lambda(\gamma) = e^{\Lambda(\gamma)}$ is an eigenvalue of $\hat{J}_\gamma$ and $\hat{J}_\gamma$ is $e^{-\Lambda(\gamma)}$-recurrent. We will need a little more since we are also interested in the $e^{-\Lambda(\gamma)}$-transient case.

LEMMA 1. *We have $R(\hat{J}_\gamma) = \sup\{u : \Psi(\gamma,u) \leq 1\}$. Also, $\hat{J}_\gamma$ is $R(\hat{J}_\gamma)$-recurrent if and only if $\Psi(\gamma,\hat{R}(J_\gamma)) = 1$.*

PROOF. Note that $\Psi(\gamma,u)$ is strictly increasing in $u$. Both results follow from

$$(3) \quad \sum_{x,m} J^m((0,\hat{y}),(x,\hat{y}))e^{\gamma x}u^m = \sum_{x,m} J_\gamma^m((0,\hat{y}),(x,\hat{y}))u^m = \frac{1}{1-\Psi(\gamma,u)}. \quad \square$$

**3. Asymptotics of $\mathcal{G}_0((0,0);(\ell,0))$ when $\hat{\mathcal{I}}$ is a nearest neighbour walk.** In Theorem 4, we give the asymptotics as $\ell \to \infty$ of the steady state $\pi(\ell,\hat{y})$ of a Markov chain in terms of the asymptotics of $\mathcal{G}((0,\hat{\sigma}),(\ell,\hat{\sigma}))$, where $\mathcal{G}$ is the potential of the associated twisted (boundary-free) Markov additive chain and $\hat{\sigma}$ is some fixed state. In general, it may be impossible to obtain these asymptotics, but in this section, we determine exact asymptotics of the potential for a Markov additive chain with a particular structure. We denote the transition kernel of this particular Markov chain by $\mathcal{I}$ and its potential matrix by $\mathcal{G}_0 := \sum_{n=0}^\infty \mathcal{I}^n$. Throughout, we use the subscript 0 to emphasize that this process may be killed when the Markovian part hits the fixed state $\hat{\sigma} = 0$.

Let $(\mathcal{V}_0,\mathcal{Z}_0)$ be a Markov additive process with transition kernel $\mathcal{I}$. We will assume that $\mathcal{I}$ is the transition kernel of a two-dimensional random walk when the Markovian component is positive and that the Markovian component is a nearest neighbor random walk on $\hat{S} = \{0,1,2,\ldots\}$ with a killing probability $\kappa \geq 0$ when the Markovian part is zero.

More precisely, we assume for $\mathcal{I}$,

$$\mathcal{I}((0,y),(x,y+z)) = \mathcal{I}((0,1),(x,1+z)) \qquad \text{for } y > 0$$

and for $\hat{\mathcal{I}}$,

$$\hat{\mathcal{I}}(y, y+1) = \begin{cases} p, & \text{for } y > 0, \\ p_0, & \text{for } y = 0, \end{cases}$$

$$\hat{\mathcal{I}}(y, y) = \begin{cases} s, & \text{for } y > 0, \\ s_0, & \text{for } y = 0, \end{cases}$$

$$\hat{\mathcal{I}}(y, y-1) = \begin{cases} q, & \text{for } y > 0, \\ 0, & \text{for } y = 0 \end{cases}$$

with $p > 0$, $q > 0$, $p + q + s = 1$, $p_0 > 0$ and $\kappa = 1 - p_0 - s_0 \geq 0$, where $\kappa$ is the probability that the process is killed when the Markovian part is zero.

To complete the specification of $\mathcal{I}$, assume that we are given the following transforms:

$$P(z) = \sum_x \mathcal{I}((0, y), (x, y+1)) z^x,$$

$$S(z) = \sum_x \mathcal{I}((0, y), (x, y)) z^x,$$

(4) $\qquad Q(z) = \sum_x \mathcal{I}((0, y), (x, y-1)) z^x \qquad \text{for } y > 0,$

$$P_0(z) = \sum_x \mathcal{I}((0, 0), (x, 1)) z^x,$$

$$S_0(z) = \sum_x \mathcal{I}((0, 0), (x, 0)) z^x \qquad \text{for } y = 0,$$

which are assumed to be finite in some neighborhood of 1. The first three transforms describe the behavior above, and the last two on the axis. Note that $P(1) = p$, $S(1) = s$, $Q(1) = q$, $P_0(1) = p_0$, $S_0(1) = s_0$ and that the horizontal drift above the $x$-axis is given by

(5) $\qquad\qquad\qquad d_+ = (P'(1) + Q'(1) + S'(1)).$

Also note that the five functions $P(z)$, $Q(z)$, $S(z)$, $P_0(z)$ and $S_0(z)$ are analytic in some annulus $D^{r,R}$, where $D^{r,R} = \{z \in \mathbb{C} | r < |z| < R\}$, $\mathbb{C}$ is the complex plane, $r < 1 < R$, and the right-hand sides form the Laurent series representations of the functions on the annulus $D^{r,R}$.

We will be particularly interested when $\hat{\mathcal{I}}$ is either 1-transient and 1-null recurrent, which corresponds to the bridge case. Under these conditions, Proposition 1 gives the exact asymptotics of $\mathcal{G}_0((0,0); (\ell, 0))$, which is simply the expected number of visits to some distant point $(\ell, 0)$ given that $(\mathcal{V}_0, \mathcal{Z}_0)$ started at the origin.

Recall that a discrete (substochastic) density $h$ on the integers has period $r \geq 1$ if $r = gcd\{u : h(u) > 0\}$, where $gcd$ denotes the greatest common divisor of the set. Equivalently, $r$ is the largest integer such that the support of $h$ is contained in the set $\{kr, \text{ where } k \in \mathbb{Z}\}$. We will say that $H(z) \equiv \sum_{\{x \in \mathbb{Z}\}} h(x) z^x$ has period $r$ if $h$ has period $r$. If $H$ has period $r$, then $H(\exp(2k\pi i/r)) = H(1)$ for all integers $k$, but $|H(\cdot)|$ is strictly less than $H(1)$ on the rest of the unit circle.

PROPOSITION 1. *If $p = q$, then the spectral radius of $\hat{\mathcal{I}}$ is one. In addition, suppose that $d_+ > 0$ and that the following aperiodicity condition holds. If $s > 0$, the period $r_{U+D}$ of the transform $P(z)Q(z)/(pq)$ and the period $r_S$ of the transform $S(z)/s$ must be relatively prime; otherwise, if $s = 0$, then $r_{U+D}$ must be one. Under these three conditions,*

$$\mathcal{G}_0((0,0);(\ell,0)) := \sum_{n=0}^{\infty} \mathcal{I}^n((0,0);(\ell,0)) \sim \begin{cases} C_+ \ell^{-3/2}, & \text{for } \kappa > 0, \\ C_0 \ell^{-1/2}, & \text{for } \kappa = 0, \end{cases}$$

*where*

$$C_+ = \frac{p_0}{\kappa^2} \sqrt{\frac{d_+}{2\pi(1-s)}} \quad \text{and} \quad C_0 = \sqrt{\frac{1-s}{2\pi d_+}}.$$

Before proving Proposition 1, we derive some expressions that will be used in the proof (and hold without the proposition's hypotheses). Roughly, our approach is to derive the generating function $\mathcal{G}_0(z) = \sum_{n \geq 0} \mathcal{G}_0((0,0);(n,0)) z^n$ and then extract the coefficients asymptotically. The generating function $\mathcal{G}_0(z)$ can be written as $1/(1 - F(z))$, where $F(z) = \mathrm{E}_{(0,0)}[z^{\mathcal{V}_0[\tau]}]$ and $\tau > 0$ is the number of steps until the Markov chain $\mathcal{Z}_0$ returns to zero. Note that $F(1) = 1 - \kappa$. Rather than deriving $F(z)$ directly, we derive $\mathrm{E}_{(0,0)}[u^\tau z^{\mathcal{V}_0[\tau]}]$, which slightly generalizes results in Section 6 of [20].

Temporarily, assume $s > 0$ and $s_0 > 0$. Let $U$ be the number of upward steps taken by time $\tau$, and let $f_n = P(U = n)$. Hence, $f_0 = s_0$ and

(6) $$\mathrm{E}_{(0,0)}[u^\tau z^{\mathcal{V}_0[\tau]} | U = 0] = u S_0(z)/s_0.$$

To calculate $f_n$ when $n > 0$, note that the first step must be up, which has probability $p_0$. If we ignore all steps where $\mathcal{Z}_0$ stays put, then the probability of an upward step is $p/(p+q)$, while the probability of a downward step is $q/(p+q)$. The total number of paths is $c_{n-1}$, the $(n-1)$st Catalan number; that is, $c_n$ is the number of paths of $2n$ steps on the positive $y$-axis that start and end at the origin and contain only nearest neighbor steps. Hence,

$$f_n = c_{n-1} p_0 \left(\frac{p}{p+q}\right)^{n-1} \left(\frac{q}{p+q}\right)^n.$$

It is well known that $c_n = \frac{1}{n+1}\binom{2n}{n}$ and that the series $\sum_{n=1}^{\infty} c_{n-1} u^n = (1 - \sqrt{1-4u})/2$ with a radius of convergence $\frac{1}{4}$. The square root is defined using a branch cut along the negative real axis.

Given $U = n > 0$, the first step must be up, but after that there are $2n - 1$ different places where $\mathcal{Z}_0$ can stay put for a geometric number of transitions. Hence, $\tau = 2n + \sum_{k=1}^{2n-1} E_k$, where $E_k$ are independent geometric random variables with law $P(E_k = m) = s^m(1 - s)$, $m = 0, 1, \ldots$. Each $E_k$ represents the number of times $\mathcal{Z}_0$ stays put between jumps. The transform of the associated displacement when $\mathcal{Z}_0$ stays put is $S(z)/s$. If we let $X[k]$ be the additive displacement during the $k$th transition when $\mathcal{Z}_0$ stays put,

$$\mathrm{E}[u^{E_k} z^{\sum_{k=1}^{E_k} X[k]}] = \sum_{m=0}^{\infty} s^m(1-s)u^m \left(\frac{S(z)}{s}\right)^m = \frac{1-s}{1-uS(z)},$$

provided $|uS(z)| < 1$. Hence,

$$\mathrm{E}_{(0,0)}[u^\tau z^{\mathcal{V}_0[\tau]} | U = n > 0] \tag{7}$$
$$= \frac{P_0(z)}{p_0} \left(\frac{P(z)}{p}\right)^{n-1} \left(\frac{Q(z)}{q}\right)^n u^{2n} \left(\frac{1-s}{1-uS(z)}\right)^{2n-1},$$

assuming $|uS(z)| < 1$.

It follows that if $|uS(z)| < 1$,

$$\mathrm{E}_{(0,0)}[u^\tau z^{\mathcal{V}_0[\tau]}]$$

$$= f_0 u \frac{S_0(z)}{s_0} + \sum_{n=1}^{\infty} f_n \frac{P_0(z)}{p_0} \left(\frac{P(z)}{p}\right)^{n-1} \left(\frac{Q(z)}{q}\right)^n u^{2n} \left(\frac{1-s}{1-uS(z)}\right)^{2n-1}$$

$$= uS_0(z) + \sum_{n=1}^{\infty} c_{n-1} p_0 \left(\frac{p}{p+q}\right)^{n-1} \left(\frac{q}{p+q}\right)^n$$
$$\times \frac{P_0(z)}{p_0} \left(\frac{P(z)}{p}\right)^{n-1} \left(\frac{Q(z)}{q}\right)^n u^{2n} \left(\frac{1-s}{1-uS(z)}\right)^{2n-1}$$

$$= uS_0(z) + \frac{1-uS(z)}{1-s} \frac{p+q}{p} \frac{p}{P(z)} P_0(z)$$
$$\times \sum_{n=1}^{\infty} c_{n-1} \left(\frac{P(z)}{p+q} \frac{Q(z)}{p+q}\right)^n u^{2n} \left(\frac{1-s}{1-uS(z)}\right)^{2n}$$

$$= uS_0(z) + \frac{1-uS(z)}{1-s} \frac{p+q}{P(z)} P_0(z) \tag{8}$$
$$\times \frac{1}{2}\left(1 - \sqrt{1 - 4\frac{P(z)}{p+q}\frac{Q(z)}{p+q} u^2 \left(\frac{1-s}{1-uS(z)}\right)^2}\right),$$

where the last summation required the term inside the square root be positive. Note that the last three expressions also hold when $s = 0$ or $s_0 = 0$ so we can resume assuming that $s \geq 0$ and $s_0 \geq 0$.

Before using (8) to compute $F(z)$, we will use it to compute the radius of convergence of $\mathrm{E}_{(0,0)}[u^\tau]$, which will be used in the proof of Proposition 1. Let $z = 1$ in (8), and notice that $\mathrm{E}_{(0,0)}[u^\tau]$ exists iff

(9) $$|us| < 1$$

and

(10) $$\left| \frac{p}{p+q} \frac{q}{p+q} u^2 \left( \frac{1-s}{1-us} \right)^2 \right| \leq \frac{1}{4},$$

in which case

(11) $$\mathrm{E}_{(0,0)}[u^\tau] = s_0 u + \frac{(1-us)}{(1-s)} \frac{p+q}{p} p_0$$
$$\times \frac{1}{2} \left( 1 - \sqrt{1 - 4 \frac{p}{p+q} \frac{q}{p+q} u^2 \left( \frac{1-s}{1-us} \right)^2} \right).$$

It is easy to see there exists positive $u < 1/s$ such that (10) is not satisfied. Since $\mathrm{E}_{(0,0)}[u^\tau]$ is an increasing function of $u$, there is a unique, positive $R < 1/s$ such that

(12) $$\frac{p}{p+q} \frac{q}{p+q} R^2 \left( \frac{1-s}{1-Rs} \right)^2 = \frac{1}{4},$$

and $R$ must be the radius of convergence of the power series in $u$ given by $\mathrm{E}_{(0,0)}[u^\tau]$.

To obtain F(z), let $u = 1$ in (8), which gives

(13) $$F(z) = S_0(z) + \frac{1 - S(z)}{1 - s} \frac{p+q}{P(z)} \frac{P_0(z)}{2}$$
$$\times \left( 1 - \sqrt{1 - 4 \frac{P(z)}{p+q} \frac{Q(z)}{p+q} \left( \frac{1-s}{1-S(z)} \right)^2} \right)$$

(14) $$= A(z) - B(z) \sqrt{C(z)} \sqrt{1-z},$$

where

$$A(z) = S_0(z) + B(z),$$
$$B(z) = \frac{1 - S(z)}{1 - s} \frac{p+q}{P(z)} \frac{P_0(z)}{2},$$
$$C(z) = \frac{1}{1-z} \left( 1 - 4 \frac{P(z)}{p+q} \frac{Q(z)}{p+q} \left( \frac{1-s}{1-S(z)} \right)^2 \right).$$

For future reference, note that

$$\mathrm{E}_{(0,y)}[z^{\mathcal{V}_0[\tau]}] = \left[ \frac{1-S(z)}{1-s} \frac{p+q}{P(z)} \frac{1}{2} \right.$$
$$\left. \times \left(1 - \sqrt{1 - 4\frac{P(z)}{p+q}\frac{Q(z)}{p+q}\left(\frac{1-s}{1-S(z)}\right)^2}\right) \right]^y \tag{15}$$

for $y = 1, 2, \ldots,$

which can be obtained by using (13), $F(z) = S_0(z) + P_0(z)\mathrm{E}_{(0,1)}[z^{\mathcal{V}_0[\tau]}]$ and $\mathrm{E}_{(0,y)}[z^{\mathcal{V}_0[\tau]}] = \mathrm{E}_{(0,1)}[z^{\mathcal{V}_0[\tau]}]^y$.

We will use the following lemma in the proof of Proposition 1.

LEMMA 2. *If $p = q$, $d_+ > 0$, and the aperiodicity condition given in Proposition 1 holds, then $\sqrt{C(z)}$ is analytic in an annulus $D^{r,R}$, where $r < 1 < R$.*

PROOF. First, we show that $\sqrt{C(z)}$ is analytic for $z \in D^{r,R} \setminus \{1\}$. On the contrary, suppose $\sqrt{C(z)}$ has a singularity on the unit circle at $z_0 \neq 1$. This means

$$4\frac{P(z_0)}{p+q}\frac{Q(z_0)}{p+q}\left(\frac{1-s}{1-S(z_0)}\right)^2 = 1. \tag{16}$$

Since $|P(z_0)| \leq p$, $|Q(z_0)| \leq q$ and $|1 - S(z_0)|^2 \geq (1 - |S(z_0)|)^2 \geq (1-s)^2$, the only way that (16) can hold is if $|P(z_0)| = p$, $|Q(z_0)| = q$ and $|1 - S(z_0)| = 1 - s$. The last equality holds only if $S(z_0) = s$. Moreover, $P(z_0)Q(z_0)/(pq) = 1$, which means that $P(z)Q(z)/(pq)$ must be periodic with a period $r_{U+D} > 1$ since $z_0 \neq 1$. Similarly, if $s > 0$, then $S(z)/s$ must be periodic with a period $r_S > 1$.

Suppose that $s > 0$. By the periodicity of $S(z)/s$, there exists a positive integer $k_S < r_S$ so that $\theta = 2\pi\frac{k_S}{r_S}$ and $z_0 = \exp(\theta i)$. Similarly, there exists a positive $k_{U+D} < r_{U+D}$ so that $\theta = 2\pi\frac{k_{U+D}}{r_{U+D}}$. Thus,

$$0 < \frac{k_S}{r_S} = \frac{k_{U+D}}{r_{U+D}} = f < 1.$$

As a consequence of Euclid's algorithm, there are integers $m_S$, $m_{U+D}$ such that $m_S \cdot r_S + m_{U+D} \cdot r_{U+D} = gcd\{r_S, r_{U+D}\} = 1$, where the last equality follows from $r_S$ and $r_{U+D}$ being relatively prime. However, this means $m_S \cdot k_S + m_{U+D} \times k_{U+D} = f$, which is a contradiction since the left-hand side is an integer, but $0 < f < 1$. If $s = 0$, then we come to the same conclusion because $r_{U+D}$ has period one by hypothesis. Hence, $\sqrt{C(z)}$ is free of singularities on $z \in D^{r,R} \setminus \{1\}$; to complete the proof, we need only show that $\sqrt{C(z)}$ is also analytic at $z = 1$.

By l'Hôpital's rule,

$$\lim_{z \to 1} C(z) = 4\left(\frac{P'(1)}{p+q}\frac{q}{p+q} + \frac{p}{p+q}\frac{Q'(1)}{p+q} + \frac{p}{p+q}\frac{q}{p+q}2S'(1)\left(\frac{1}{(1-s)}\right)\right)$$

$$= 2\left(\frac{P'(1)}{p+q} + \frac{Q'(1)}{p+q} + \frac{S'(1)}{(1-s)}\right) \quad \text{since } p = q$$

(17) $$\equiv 2d_+/(1-s).$$

By hypothesis, $d_+ > 0$ so the limit exists as $z \to 1$ and is in the domain of the square root. Hence, it follows that $\sqrt{C(z)}$ is analytic in some annulus $D^{r,R}$. □

PROOF OF PROPOSITION 1. Since $p = q$, $R = 1$ satisfies (12) and obviously $R < 1/s$. So the radius of convergence of $\mathrm{E}_{(0,0)}(u^\tau)$ is 1, as desired. Note also that by (11), $\mathrm{E}_{(0,0)}(1^\tau) = \Pr\{\tau < \infty\} = s_0 + p_0 = 1 - \kappa$,

$$\mathcal{G}_0(z) = \sum_{n=0}^{\infty} \mathcal{G}_0((0,0);(n,0))z^n = \frac{1}{1 - F(z)}$$

$$= (1 - A(z) + B(z)\sqrt{C(z)}\sqrt{1-z})^{-1}$$

$$= \frac{1 - A(z) - B(z)\sqrt{C(z)}\sqrt{1-z}}{(1 - A(z))^2 - B^2(z)C(z)(1-z)}$$

$$= r(z) + \frac{r_2(z)}{r_0(z)}\sqrt{C(z)}\sqrt{1-z},$$

where

$$r_0(z) = (1 - A(z))^2 - B^2(z)C(z)(1-z),$$

$$r_1(z) = 1 - A(z), \qquad r_2(z) = -B(z), \qquad r(z) = \frac{r_1(z)}{r_0(z)}.$$

$S_0(z)$, $P_0(z)$, $S(z)$, $P(z)$, $Q(z)$ are analytic in some annulus $D^{r,R}$ as is $\sqrt{C(z)}$ by Lemma 2. Furthermore,

$$r_0(z) = (1 - S_0(z) - B(z))^2 - B^2(z)\left(1 - 4\frac{P(z)}{p+q}\frac{Q(z)}{p+q}\left(\frac{1-s}{1-S(z)}\right)^2\right)$$

$$= (1 - S_0(z))^2 - 2B(z)(1 - S_0(z)) + 4B^2(z)\frac{P(z)}{p+q}\frac{Q(z)}{p+q}\left(\frac{1-s}{1-S(z)}\right)^2$$

$$= (1 - S_0(z))^2 - \frac{1 - S(z)}{1 - s}\frac{p+q}{P(z)}P_0(z)(1 - S_0(z))$$

$$+ 4\left(\frac{1}{2}\frac{1 - S(z)}{1 - s}\frac{p+q}{P(z)}P_0(z)\right)^2\frac{P(z)}{p+q}\frac{Q(z)}{p+q}\left(\frac{1-s}{1-S(z)}\right)^2$$

$$= (1 - S_0(z))^2 - \frac{1 - S(z)}{1 - s} \frac{p + q}{P(z)} P_0(z)(1 - S_0(z)) + P_0(z)^2 \frac{Q(z)}{P(z)}$$

$$\equiv \frac{D(z)}{P(z)},$$

where $D(z) = P(z)(1 - S_0(z))^2 - \frac{1-S(z)}{1-s}(p+q)P_0(z)(1 - S_0(z)) + P_0(z)^2 Q(z)$.

Since

$$r(z) = \frac{P(z)}{D(z)}\left(1 - S_0(z) - \frac{1}{2}\left(\frac{1 - S(z)}{1 - s}\right)\frac{p + q}{P(z)} P_0(z)\right),$$

it follows that $r(z) = r_3(z)/D(z)$, where $r_3(z)$ is analytic on $D^{r,R}$. Similarly,

$$\frac{r_2(z)}{r_0(z)} = \frac{P(z)}{D(z)}\left(-\frac{1}{2}\left(\frac{1 - S(z)}{1 - s}\right)\frac{p + q}{P(z)} P_0(z)\right)$$

so $r_2(z)/r_0(z) = r_4(z)/D(z)$, where $r_4(z)$ is analytic on $D^{r,R}$.

CASE 1. Assume $\kappa > 0$. Note that $D(z)$ has no zeroes in $\{r < |z| \leq 1\}$, where $r < 1$. Suppose it did! It would follow from the above that $1/(1 - F(z))$ would have a pole inside $\{r < |z| \leq 1\}$. This is impossible because $F(1) < 1$ [so $|F(z)| < 1$ if $|z| \leq 1$] when the killing probability $\kappa$ is greater than zero. Since $D(z)$ only has a finite number of zeroes in a compact set, it follows that there exists an annulus $D^{r,R_D}$ with $R_D > 1 > r$, where $1/D(z)$ is analytic. Hence,

$$\mathcal{G}_0((0,0);(\ell,0)) = \int_\gamma \frac{\mathcal{G}_0(z)}{z^{\ell+1}}\, dz,$$

where $\gamma$ is any positively oriented circle that encloses zero inside the domain of convergence of $\mathcal{G}_0(z)$. Since $r(z)$ is analytic on the punctured disk of radius $R_D > 1$,

$$\int_\gamma \frac{r(z)}{z^{\ell+1}}\, dz = \int_{\gamma^+} \frac{r(z)}{z^{\ell+1}}\, dz \to 0 \quad \text{exponentially fast,}$$

where $\gamma^+$ is a circle of radius greater than one. Consequently, the coefficient of $z^\ell$ in the Laurent expansion of $r(z)$ decays exponentially fast as $\ell \to \infty$ and can be neglected.

Moreover,

$$\frac{r_2(z)}{r_0(z)}\sqrt{C(z)}\sqrt{(1-z)} = -C\sqrt{(1-z)} + \mathcal{O}(|1-z|)$$

as $z \to 1$ in $D_{\alpha,\delta}$, where $D_{\alpha,\delta} = \{z \in \mathbb{C} : |z| \leq 1 + \delta, |\arg(z-1)| \geq \alpha\}$, where $\delta > 0$, $1 + \delta < R$, $0 < \alpha < \pi/2$ and where

$$C = -\sqrt{C(1)}r_2(1)/r_0(1) = \frac{p_0}{(1-p_0-s_0)^2}\sqrt{\frac{2d_+}{1-s}} = \frac{p_0}{\kappa^2}\sqrt{\frac{2d_+}{1-s}}.$$

If $1 + \delta < R_D$, then $\frac{r_2(z)}{r_0(z)}\sqrt{C(z)} = \frac{r_4(z)}{D(z)}\sqrt{C(z)}$ is analytic in $D_{\alpha,\delta}$. By Theorem 16.8 in [24],

$$\mathcal{G}_0((0,0);(\ell,0)) \sim -C\binom{\ell - 3/2}{\ell}$$
$$= -C\Gamma[\ell + 1 - 3/2]/(\Gamma[\ell + 1]\Gamma[1 - 3/2])$$
$$\sim -C\Gamma[\ell - 1/2]/(\Gamma[\ell + 1]\Gamma[-1/2]).$$

Recall $\sqrt{\pi} = \Gamma[1/2] = (-1/2)\Gamma[-1/2]$ and the fact that $\Gamma[\ell - 1/2]/\Gamma[\ell + 1] \sim \ell^{-3/2}$ as $\ell \to \infty$ (see (5.02) in [19]). This gives

$$\mathcal{G}_0((0,0);(\ell,0)) \sim \frac{C}{2\sqrt{\pi}}\ell^{-3/2} = \frac{p_0}{\kappa^2}\sqrt{\frac{d_+}{2\pi(1-s)}}\ell^{-3/2}.$$

CASE 2. Assume $\kappa = 0$. Note that in this null recurrent case, $D(1) = 0$. Calculation shows $D'(1) = p_0^2 d_+$ which is positive by hypothesis so $r_2(z)/r_0(z)$ has a simple pole at $z = 1$. Consequently, the expansion of $r_2(z)/r_0(z)$ around 1 to first order is

$$\frac{p}{p_0^2 d_+(z-1)}\left(-\frac{1}{2}\left(\frac{1-S(1)}{1-s}\right)\frac{p+q}{p}P_0(1)\right) = -\frac{1-s}{2p_0 d_+(z-1)}.$$

To show that $r(z)$ is again negligible, note that $r_1(1) = 0$, which cancels the zero of order one in $D(z)$ at $z = 1$. Hence, $r(z)$ is analytic at 1 so this term may be neglected as in the previous case. Therefore, in the neighborhood of 1,

$$\mathcal{G}_0(z) \sim \frac{-(1-s)}{2p_0 d_+(z-1)}\sqrt{C(z)}(1-z)^{1/2}$$
$$\sim \frac{1-s}{2p_0 d_+}\sqrt{C(1)}(1-z)^{-1/2} = c(1-z)^{-1/2},$$

where $c = (1-s)/(p_0\sqrt{2d_+})$.

Again by Theorem 16.8 in [24],

$$\mathcal{G}_0((0,0);(\ell,0)) \sim c\binom{\ell - 1/2}{\ell}$$
$$\sim c\Gamma[\ell + 1/2]/(\Gamma[\ell + 1]\Gamma[1/2]).$$

Recall $\sqrt{\pi} = \Gamma[1/2]$ and the fact that $\Gamma[\ell + 1/2]/\Gamma[\ell + 1] \sim \ell^{-1/2}$ as $\ell \to \infty$ (see (5.02) in [19]). This gives

$$\mathcal{G}_0((0,0);(\ell,0)) \sim \frac{c}{\sqrt{\pi}}\ell^{-1/2} = \frac{1-s}{p_0\sqrt{2\pi d_+}}\ell^{-1/2}.$$

□

The next proposition gives conditions for $F(z)$ to be aperiodic. Note that the aperiodicity condition given in Proposition 1 is more stringent than the condition given in the following proposition.

PROPOSITION 2. *If the periods of $S_0(z)$ (if $s_0 > 0$), $S(z)$ (if $s > 0$), $P_0(z)Q(z)/(p_0q)$ and $P(z)Q(z)/(pq)$ are relatively prime, then $F(z)$ has period one.*

PROOF. Consider the case when $s_0 > 0$ and $s > 0$. Define $G(z, n) = \mathrm{E}_{(0,0)}[z^{\mathcal{V}_0[\tau]}|U = n]$. Let $|z_0| = 1$. To have $F(z_0) = F(1)$, we must have $G(z_0, n) = G(1, n)$ for all choices of $n$. For $n = 0$, setting $u = 1$ in (6) implies that $S_0(z_0)/s_0 = 1$. For $n = 1$, setting $u = 1$ in (7) implies that $P_0(z_0)Q(z_0)/(p_0q) = 1$ and $S(z_0)/s = 1$. Similarly, for $n = 2$, we need $P_0(z_0)P(z_0)Q(z_0)^2/(p_0pq^2) = 1$ and $S(z_0) = 1$. Thus, we must have $S_0(z_0)/s_0 = 1$, $S(z_0)/s = 1$, $P_0(z_0)Q(z_0)/(p_0q) = 1$ and $P(z_0)Q(z_0)/(pq) = 1$. By an argument similar to the proof of Lemma 2, we must have $z_0 = 1$. Consequently, $F(z)$ is aperiodic. □

**4. Spectral radius of $\hat{J}_\gamma$ for a nearest neighbor random walk.** As in the beginning of Section 2, let $\hat{J}_\gamma$ be the Feynman–Kac transform associated with a Markov additive process with transition kernel $J$. Under certain assumptions on $J$, Proposition 3 below gives the spectral radius of $\hat{J}_\gamma$. Proposition 3 is Proposition 10 in [10]. Instead of the proof in [10], which appealed to the theory of large deviations, we give an algebraic proof using the results in Section 3 and the generalized Ney–Nummelin representation for $r(\hat{J}_\gamma)$ in Lemma 1.

The assumptions needed on $J$ are identical to the assumptions on $\mathcal{I}$ in the beginning of Section 3, except that there is no killing; that is, $\kappa = 0$. In the proof of Proposition 3, we crudely twist the kernel $J$ obtaining a kernel denoted by $\mathcal{I}$ that may have a killing on the boundary. After giving the proof of Proposition 3, we will show how the crude twist can be refined to give a harmonic function for $J$ for the nearest neighbor case. This approach will be generalized in the next section.

As in [10], define

$$
\begin{aligned}
R^+(\gamma, \beta) &= \sum_{x', y'} J((x, y); (x', y'))e^{\gamma(x'-x)}e^{\beta(y'-y)} \qquad y > 0 \quad \text{and} \\
R^-(\gamma, \beta) &= \sum_{x', y'} J((x, 0); (x', y'))e^{\gamma(x'-x)}e^{\beta y'}.
\end{aligned}
$$
(18)

Note that $R^+ - 1$ and $R^- - 1$ are identical to the functions $M^+$ and $M^-$ used in our companion paper [9] and that $R^+$ is strictly convex.

PROPOSITION 3 (Ignatiouk-Robert). *For each real $\gamma$, let $\beta_0 = \beta_0(\gamma)$ be the unique value of $\beta$ that minimizes $R^+(\gamma, \beta)$.*

1. *If the bridge condition $R^-(\gamma, \beta_0) \leq R^+(\gamma, \beta_0)$ holds, then $r(\hat{J}_\gamma) = R^+(\gamma, \beta_0)$.*
2. *If the jitter condition $R^-(\gamma, \beta_0) > R^+(\gamma, \beta_0)$ holds, then there is a unique $\beta_1(\gamma) \leq \beta_0(\gamma)$ such that $R^-(\gamma, \beta_1) = R^+(\gamma, \beta_1)$ and $r(\hat{J}_\gamma) = R^+(\gamma, \beta_1)$.*

PROOF. For any $y > 0$, let $u = u(\gamma) = \hat{J}_\gamma(y, y+1)$, $d = d(\gamma) = \hat{J}_\gamma(y, y-1)$ and $s = s(\gamma) = \hat{J}_\gamma(y, y)$. Hence, $R^+(\gamma, \beta) = ue^\beta + s + de^{-\beta}$, $\exp(\beta_0) = \sqrt{d/u}$ and $R^+(\gamma, \beta_0) = s + 2\sqrt{ud}$.

Assume the bridge condition holds. Define

$$\mathcal{I}((x,y);(x',y')) = \frac{J((x,y);(x',y'))e^{\gamma(x'-x)}e^{\beta_0(y'-y)}}{f},$$

where $f = R^+(\gamma, \beta_0) = s + 2\sqrt{ud}$ is chosen so that $\sum_{(x',y')} \mathcal{I}((x,y);(x',y')) = 1$ for all $y > 0$. Note that $\mathcal{I}$ is the probability transition kernel of a Markov additive process $(\mathcal{V}_0, \mathcal{Z}_0)$ with a possible killing at $y = 0$ as described in the last section [since $R^-(\gamma, \beta_0) \leq R^+(\gamma, \beta_0)$]. Also note that the choice of $\beta_0$ forces $p = q$.

By Lemma 1, $R(\hat{J}_\gamma) = \sup\{u : \Psi(\gamma, u) \leq 1\}$. Next note $\Psi(\gamma, u) = \mathrm{E}_{(0,0)}[(uf)^\tau]$ (where $\tau$ is the first time $\mathcal{Z}_0$ returns to zero) since all the factors involving $\exp(\beta_0)$ cancel out over trajectories which start and return to $y = 0$. Therefore, for a given $\gamma$,

$$R(\hat{J}_\gamma) = \sup\{u : \mathrm{E}_{(0,0)}[(uf)^\tau] \leq 1\}.$$

However with $p = q$, using (12), we see $R(\hat{J}_\gamma) = 1/f$; that is, the spectral radius of $\hat{J}_\gamma$ is $R^+(\gamma, \beta_0)$.

Now assume the jitter condition holds. By the argument following Proposition 9 in [10], we see there must exist a unique $\beta_1(\gamma) < \beta_0(\gamma)$ such that $R^-(\gamma, \beta_1) = R^+(\gamma, \beta_1)$. Moreover, $R^+(\gamma, \beta)$ is strictly decreasing for $\beta \leq \beta_0(\gamma)$ so the derivative of $R^+(\gamma, \beta)$ is negative at $\beta_1$.

Define

$$\mathcal{I}((x,y);(x',y')) = \frac{J((x,y);(x',y'))e^{\gamma(x'-x)}e^{\beta_1(y'-y)}}{f}$$

and $f = R^+(\gamma, \beta_1)$; thus, $\mathcal{I}$ is the probability transition kernel of a Markov additive process $(\mathcal{V}_0, \mathcal{Z}_0)$ without a killing at $y = 0$ as described in the last section. Again the spectral radius of $\hat{J}_\gamma$ is $f$ times the spectral radius of $\hat{\mathcal{I}}$. However, the mean drift of $\mathcal{Z}_0$ is obtained by taking the derivative with respect to $\beta$ of $R^+(\gamma, \beta)/f = (u(\gamma)e^\beta + s(\gamma) + d(\gamma)e^{-\beta})/f$ and evaluating at $\beta_1$. We already know the derivative of $R^+(\gamma, \beta)$ is negative at $\beta_1$ so we conclude $\mathcal{Z}_0$ is positive recurrent and has spectral radius one. We conclude that the spectral radius of $\hat{J}_\gamma$ is $f = R^+(\gamma, \beta_1)$. □

We have already remarked that the spectral radius of $\hat{J}_\gamma$ is a convex function and since the spectral radius of $\hat{J}_0 \equiv \hat{J}$ is one, it follows that there is at most one choice of $\gamma > 0$ such that $r(\hat{J}_\gamma) = 1$. We will call this point $\alpha$. Now, we show how to refine the crude twist to find a harmonic function for $J$. First, we consider the jitter case.

*The jitter case* $R^-(\alpha, \beta_0) > R^+(\alpha, \beta_0)$. In this case the function $\exp(\beta_1(\alpha)y)$ is a right eigenfunction for $\hat{J}_\alpha$ with eigenvalue $R^+(\alpha, \beta_1)$. Hence, $h(x,y) = \exp(\alpha x)\exp(\beta_1(\alpha)y)$ is harmonic for $J$. We can now perform the $h$ transform of $J$ to get the kernel $\mathcal{J}$ and we denote the $h$-transformed Markov chain by $(\mathcal{V}, \mathcal{Z})$. The Markovian component has kernel $\hat{\mathcal{J}}$, where

$$p = \hat{\mathcal{J}}(y, y+1) = u\exp(\beta_1(\alpha)),$$
$$q = \hat{\mathcal{J}}(y, y-1) = d\exp(-\beta_1(\alpha)) \quad \text{and} \quad s = \hat{\mathcal{J}}(y,y),$$
$$p_0 = \hat{\mathcal{J}}(0,1) = u_0\exp(\beta_1(\alpha)) \quad \text{and} \quad s_0 = \hat{\mathcal{J}}(0,0).$$

Note that $p_0 + s_0 = 1$ and $p \leq q$; otherwise, $\hat{\mathcal{J}}$ would not have spectral radius one.

*The bridge case* $R^-(\alpha, \beta_0) \leq R^+(\alpha, \beta_0)$. With the parameters $\alpha$ and $\beta_0(\alpha)$ we construct the kernel $\mathcal{I}$ which has spectral radius one. Now remark that the function $a_0(y) = (1 + \kappa y/p_0)$ is harmonic for $\hat{\mathcal{I}}$ so, in fact, the function $h(x,y) = \exp(\alpha x)\exp(\beta_0(\alpha)y)a_0(y)$ is harmonic for $J$. We can now perform the $h$ transform of $J$ to get the kernel $\mathcal{J}$ and we denote the $h$-transformed Markov chain by $(\mathcal{V}, \mathcal{Z})$. The Markovian component has kernel $\hat{\mathcal{J}}$, where

$$p(y) \equiv \hat{\mathcal{J}}(y, y+1) = u\frac{a_0(y+1)}{a_0(y)} = u\frac{1 + \kappa(y+1)/p_0}{1 + \kappa y/p_0},$$
$$q(y) \equiv \hat{\mathcal{J}}(y, y-1) = u\frac{a_0(y-1)}{a_0(y)} = u\frac{1 + \kappa(y-1)/p_0}{1 + \kappa y/p_0},$$
$$s(y) \equiv \hat{\mathcal{J}}(y,y) = s,$$

for $y > 0$ and

$$\hat{\mathcal{J}}(0,1) = p_0\frac{a_0(1)}{a_0(0)} = p_0(1 + \kappa/p_0), \qquad \hat{\mathcal{J}}(0,0) = s_0.$$

If the bridge condition holds with an inequality (and $\kappa > 0$), then the kernel $\hat{\mathcal{J}}$ is transient. To see this we remark that there is some probability of drifting to plus infinity without ever hitting zero because the criterion for transience is that

$$\sum_{n=1}^{\infty} \prod_{k=1}^{n} \frac{q(1)\cdots q(k)}{p(1)\cdots p(k)} < \infty.$$

By telescoping, the above sum is

$$\sum_{n=1}^{\infty} \frac{p_0(p_0+\kappa)}{(p_0+\kappa n)(p_0+\kappa(n+1))} < \infty.$$

If the bridge condition holds with equality (and $\kappa = 0$), then $\beta_0(\alpha) = \beta_1(\alpha)$ and the function $h(x,y) = \exp(\alpha x)\exp(\beta_1(\alpha)y)$ is harmonic for $J$. The $h$-transformed Markov chain $(\mathcal{V}, \mathcal{Z})$ has a Markovian component with kernel $\hat{\mathcal{J}}$ with $p = q$, which is null recurrent.

**5. The $h$-transform approach and a ratio limit theorem.** The last section gave the exact asymptotics for Markov additive processes whose Markovian part is a random walk with a boundary. In this section we allow the Markovian part to be a general Markov chain on a countable state space. We assume that $J$ has a positive harmonic function $h$; that is, $h > 0$ and $Jh = h$. Furthermore, we assume that $h$ has the form $h(z) = e^{\alpha z_1}\hat{h}(\hat{z})$, where $\alpha > 0$. We use $h$ to construct the twisted process $(\mathcal{V}[n], \mathcal{Z}[n])$ having transition kernel $\mathcal{J}(z,x) \equiv J(z,x)h(x)/h(z)$. We use caligraphic letters for the twisted process. The probabilities of the two processes for a sequence of states $x[0], x[1], \ldots, x[n]$ in $\mathbb{Z} \times \hat{S}$ are related via the following change of measure:

$$\Pr\{(V[n], Z[n]) = x[n], \ldots, (V[1], Z[1]) = x[1] | (V[0], Z[0]) = x[0]\}$$
$$= \Pr\{(\mathcal{V}[n], \mathcal{Z}[n]) = x[n], \ldots, (\mathcal{V}[1], \mathcal{Z}[1]) = x[1] | (\mathcal{V}[0], \mathcal{Z}[0]) = x[0]\}$$
$$\times h(x[0])/h(x[n]).$$

Our goal is to investigate the asymptotics of $G(z, (\ell, \hat{x}))$ as $\ell \to \infty$, where $G$ is the Green function of $J$,

$$G(z,x) \equiv \sum_{n=0}^{\infty} J^n(z,x),$$

which gives the expected number of visits to any state $x$ starting from any state $z$. Unfortunately, we cannot obtain sharp asymptotics as in Section 3. Our main result is a ratio limit theorem for $G(z, (\ell + r, \hat{x}))/G(w, (\ell, \hat{y}))$ as $\ell \to \infty$.

Define the Green's function

$$\mathcal{G}(z,x) = \sum_{n=0}^{\infty} \mathcal{J}^n(z,x)$$

of $\mathcal{J}$. The Green's function $G$ of $J$ are related by $G(z,x) = \mathcal{G}(z,x)h(z)/h(x)$. We find it easier to investigate the asymptotics of $G$ by studying $\mathcal{G}$ since $\mathcal{V}[n]$ will be assumed to drift to $+\infty$. As a consequence of the drift, $(\mathcal{V}, \mathcal{Z})$ is transient and $\mathcal{G}$ and $G$ are finite.

Note that the marginal processes $Z \equiv Z[0], Z[1], \ldots$ and $\mathcal{Z} \equiv \mathcal{Z}[0], \mathcal{Z}[1], \ldots$ are Markov chains, and we denote their transition kernels by $\hat{J}$ and $\hat{\mathcal{J}}$, respectively. We are interested in the case when $\hat{\mathcal{J}}$ has spectral radius 1 and is either 1-transient or null recurrent since this will be needed in analyzing the "bridges" of the modified Jackson network, which fall outside of the scope of [8, 14].

The papers [8, 14] consider the case when $\hat{\mathcal{J}}$ is positive recurrent; in this case, under reasonable assumptions, $\lim_{\ell \to \infty} \mathcal{G}(z, (\ell, \hat{x}))$ converges to a positive limit, which is a function of the invariant probability measure of $\hat{\mathcal{J}}$ and the speed at which $\mathcal{V}$ drifts to positive infinity. However, when $\hat{\mathcal{J}}$ is 1-transient or null recurrent, it does not have an invariant probability measure and $\lim_{\ell \to \infty} \mathcal{G}(z, (\ell, \hat{x})) = 0$. Instead, we consider the ratio

$$\frac{\mathcal{G}(z, (\ell, \hat{x}))}{\mathcal{G}(w, (\ell, \hat{y}))}$$

as $\ell$ goes to infinity.

The assumptions we need on the Markov additive processes $(V, Z)$ and $(\mathcal{V}, \mathcal{Z})$ for our main result are the following:

A0. The transition kernel $J$ has a positive harmonic function $h(z) = e^{\alpha z_1} \hat{h}(\hat{z})$, where $\alpha > 0$.
A1. There exists a $t > \alpha$ and $M$ such that $\mathrm{E}(\exp(tV[1])|Z[0], Z[1]) < M$ for all $Z[0]$ and $Z[1]$.
A2. $\hat{J}$ is irreducible.
A2.5. Define $T_{\hat{\sigma}}$ to be number of steps for $Z$ to return to some state $\hat{\sigma} \in \hat{S}$. We assume that the distribution of $P_{(0,\hat{\sigma})}(V[T_{\hat{\sigma}}] = \cdot)$ is not concentrated on a subgroup of the integers. (This is Condition (P(1)) in [14]).
A3. $\mathcal{V}[n] \to \infty$ almost surely as $n \to \infty$.
A4. $\hat{\mathcal{J}}$ has spectral radius 1.
A5. $\hat{\mathcal{J}}$ has an invariant measure $\varphi$, which is unique up to constant multiples.
A5.5. We assume that constant functions are the unique harmonic functions for $\hat{\mathcal{J}}$. (This is the strong Liouville property; see [24]. Note that A5.5 automatically holds when $\hat{\mathcal{J}}$ is null recurrent.)
A6. $\mathcal{J}$ is irreducible in the sense that the probability $p((0,\hat{x}), (0,\hat{y}))$ of going from $(0,\hat{x})$ to $(0,\hat{y})$ is positive for any $\hat{x}, \hat{y} \in \hat{S}$. Also, there exists an integer $N$ and $p > 0$ fixed such that for any $\hat{y}$, there exists an integer $m = m(\hat{y})$ such that $1 \le m \le N$ and $\mathcal{J}^m((0,\hat{y}), (1,\hat{y})) \ge p$.

(Assumption A6 implies A2; however, we keep both since we believe A2 is necessary, but A6 may be too strong.)

A7. There is some state $\hat{\sigma}$ and a function $\hat{f}$ such that uniformly in $s$,

$$\mathcal{G}((0,\hat{x}),(s,\hat{\sigma}))/\mathcal{G}((0,\hat{\sigma}),(s,\hat{\sigma})) \leq \hat{f}(\hat{x}),$$

where

$$\sum_{\hat{x}\in\hat{S}} \hat{\mathcal{J}}(\hat{z},\hat{x})\hat{f}(\hat{x}) < \infty$$

for all states $\hat{z}$. [A sufficient condition for A7 is that $(V,Z)$ has bounded jumps. To see this, let $p((0,\hat{x}),(0,\hat{\sigma}))$ be the probability of ever going from $(0,\hat{x})$ to $(t,\hat{\sigma})$. Clearly,

$$p((0,\hat{\sigma}),(s,\hat{\sigma})) \geq p((0,\hat{\sigma}),(0,\hat{x}))p((0,\hat{x}),(s,\hat{\sigma})).$$

Hence,

$$\frac{\mathcal{G}((0,\hat{x}),(s,\hat{\sigma}))}{\mathcal{G}((0,\hat{\sigma}),(s,\hat{\sigma}))} = \frac{p((0,\hat{x}),(s,\hat{\sigma}))}{p((0,\hat{\sigma}),(s,\hat{\sigma}))} \leq \frac{1}{p((0,\hat{\sigma}),(0,\hat{x}))}.$$

For A7 to hold it, therefore, suffices that the range of $\mathcal{J}((0,\hat{z}),(\cdot,\cdot))$ is finite for all $\hat{z}$.]

A7*. Condition A7 holds for the Green's function of $(-\mathcal{V}^*[n], \mathcal{Z}^*[n])$, where $(\mathcal{V}^*[n], \mathcal{Z}^*[n])$ is the time reversal of $(\mathcal{V}^*[n], \mathcal{Z}^*[n])$ with respect to $\varphi$. (The purpose of the minus sign in $(-\mathcal{V}^*[n], \mathcal{Z}^*[n])$ is simply to have the process drift to the right; that is, to also satisfy condition A3.)

We also use the following convention throughout. For real valued functions $f$ and $g$, let $f(\ell) \sim g(\ell)$ mean that $\lim_{\ell\to\infty} f(\ell)/g(\ell) = 1$.

THEOREM 1. *Under assumptions* A0–A7 *and* A7* *with* $w, \hat{x}, \hat{y}$ *and* $z$ *fixed,*

$$\frac{\mathcal{G}(z,(\ell,\hat{x}))}{\mathcal{G}(w,(\ell,\hat{y}))} \sim \frac{\varphi(\hat{x})}{\varphi(\hat{y})},$$

*which implies*

$$\frac{G(z,(\ell,\hat{x}))}{G(w,(\ell,\hat{y}))} \sim \frac{\hat{h}(\hat{x})\varphi(\hat{x})}{\hat{h}(\hat{y})\varphi(\hat{y})} \frac{\hat{h}(\hat{w})}{\hat{h}(\hat{z})}.$$

The proof of Theorem 1 is broken up into proving each of the following over the next few sections:

(19) $\qquad \mathcal{G}((0,\hat{y}),(\ell,\hat{y}))^{1/\ell} \sim 1,$

(20) $\qquad \mathcal{G}((t,\hat{z}),(\ell,\hat{y})) \sim \mathcal{G}((0,\hat{z}),(\ell,\hat{y})),$

$$\mathcal{G}((t,\hat{z}),(\ell,\hat{\sigma})) \sim \mathcal{G}((0,\hat{\sigma}),(\ell,\hat{\sigma})), \tag{21}$$

$$\mathcal{G}((t,\hat{z}),(\ell,\hat{y})) \sim \mathcal{G}((0,\hat{y}),(\ell,\hat{y})), \tag{22}$$

$$\frac{\mathcal{G}((0,\hat{y}),(\ell,\hat{x}))}{\mathcal{G}((0,\hat{y}),(\ell,\hat{y}))} \sim \frac{\varphi(\hat{x})}{\varphi(\hat{y})} \tag{23}$$

for all $s, \hat{z}, \hat{y}, \hat{x}$ but $\hat{\sigma}$ satisfying A7 and A7*.

Equations (22) and (23) combine to give the first result in Theorem 1, and the second result immediately follows from a change of measure. Equation (19) implies that terms of the form $ce^{-s\ell}$ make an asymptotically negligible contribution to $\mathcal{G}((0,\hat{y}),(\ell,\hat{y}))$. Using this, we derive (20), which shows that asymptotics are unaffected by changing $\mathcal{V}[0]$. To show that the asymptotics are also unaffected by changing $\mathcal{Z}[0]$, we show it for the target $(\ell,\hat{\sigma})$, where $\hat{\sigma}$ satisfies A7. Once there exists $\hat{\sigma}$ such that (21) holds (whether A7 holds or not), the result is extended to all $\hat{y} \in \hat{S}$, giving (22), which implies that the asymptotics are not affected by the starting state. To obtain (23), the same arguments are repeated on the process $(-\mathcal{V}^*[n], \mathcal{Z}^*[n])$, where * denotes the time reversal with respect to the invariant measure $\varphi$.

Since (20) is equivalent to $\mathcal{G}((0,\hat{z}),(\ell,\hat{y})) \sim \mathcal{G}((0,\hat{z}),(\ell+1,\hat{y}))$, the result might initially appear to follow directly from renewal theory using A2.5. However, consider the following: Chung ([2], page 50) gives an example of an irreducible, aperiodic Markov chain with transition matrix $P$ such that $\limsup_{n\to\infty} P^{n+1}(0,0)/P^n(0,0) = \infty$. Let $\mathcal{Z}[t]$ be that Markov chain, and consider the Markov additive process $(\mathcal{V}[t], \mathcal{Z}[t]) = (t, \mathcal{Z}[t])$. Since $\mathcal{G}((0,0),(\ell,0)) = P^\ell(0,0)$, we have $\limsup_{\ell\to\infty} \mathcal{G}((0,0),(\ell+1,0))/\mathcal{G}((0,0),(\ell,0)) = \infty$. Now modify the transition kernel of the Markov additive process $(t, \mathcal{Z}[t])$ to allow the process to either remain in the same state with probability $1/2$ or to jump as before with probability $1/2$. Since this simply doubles the expected number of visits to any state, we still have $\limsup_{\ell\to\infty} \mathcal{G}((0,0),(\ell+1,0))/\mathcal{G}((0,0),(\ell,0)) = \infty$, which violates (20). Note that A6 fails for this example.

Kesten in [11] studies ratio limit theorems for a transition kernel $P$ having spectral radius $r$. In Theorem 1 of [11], (1.4u) implies that there is a unique (up to constant multiples) positive harmonic function $h$; that is, a unique $h$ such that $Ph = rh$, while (1.4d) implies there is a unique invariant measure $\varphi$ (up to constant multiples); that is, $\varphi P = r\varphi$ (though only for the the state space specified in [11]). In this case Theorem 2 in [11] becomes

$$\lim_{\ell\to\infty} \frac{P^{\ell+k}(\hat{z},\hat{y})}{P^\ell(\hat{w},\hat{x})} = r^k \frac{h(\hat{z})\varphi(\hat{y})}{h(\hat{w})\varphi(\hat{x})}. \tag{24}$$

His ratio limit theorems with $r = 1$ are closely related to our study of the asymptotics of Markov additive processes $(\mathcal{V}, \mathcal{Z})$ satisfying A0–A7. The level crossing process $(t, \mathcal{U}[t])$ associated with a nearest neighbor Markov additive process $(\mathcal{V}, \mathcal{Z})$ is a Markov chain with kernel $P$ and spectral radius 1. Clearly,

$$(25) \quad \frac{\sum_{\hat{v}} P^{\ell}(\hat{z}, \hat{v}) \cdot \mathcal{G}((0, \hat{v}), (0, \hat{y}))}{\sum_{\hat{v}} P^{\ell}(\hat{w}, \hat{v}) \cdot \mathcal{G}((0, \hat{v}), (0, \hat{x}))} = \frac{\mathcal{G}((0, \hat{z}), (\ell, \hat{y}))}{\mathcal{G}((0, \hat{w}), (\ell, \hat{x}))},$$

so the ratio of the Green's functions will converge if the ratio limit theorem for the level crossing process holds.

Our Lemma 5 implies Lemma 4 in [11], and our assumption A6 is equivalent to (1.5) in [11] (and basically the proofs are the same). Theorem 2 in [11] requires (1.4u) and (1.4d), which essentially mean that the range of transitions are bounded—an assumption that we need to avoid because $P$ has unbounded jumps if the underlying Markov additive process has negative increments. Instead of (1.4u) and (1.4d), we simply assumed in A5 that the invariant measure is unique up to multiplication by a constant and in A5.5 that the only positive harmonic functions are the constants. However, in order to push through the argument giving (21) and its analog for the reversed process without the bounded jump assumption, we added the uniform integrability assumptions A7 and A7*. These assumptions may be too strong. However, some additional assumption is needed since Section 5.5 contains an example satisfying A0–A6, yet (21) fails.

5.1. *Proof of* (19). By the convexity of $r(\hat{J}_\gamma)$, it follows that there can be at most one point $\gamma > 0$ such that $r(\hat{J}_\gamma) = 1$; that is, when $\gamma = \alpha$. Note that $J_\alpha((x_1, \hat{y}), (y_1, \hat{y})) = \mathcal{J}((x_1, \hat{y}), (y_1, \hat{y}))$, where $\alpha$ was given in assumption A0. Consequently, $\hat{J}_\alpha^n(\hat{y}, \hat{y}) = \hat{\mathcal{J}}^n(\hat{y}, \hat{y})$. Thus, they have the same radius of convergence, that is, $R(\hat{J}_\alpha) = R(\hat{\mathcal{J}})$. By A4, the spectral radius of $\hat{\mathcal{J}} = 1$ so

$$r(\hat{J}_\alpha) = r(\hat{\mathcal{J}}) = 1; \qquad \text{hence, } R(\hat{J}_\alpha) = R(\hat{\mathcal{J}}),$$

and both kernels are either 1-recurrent, 1-null recurrent or 1-transient.

Note that

$$\Psi(\gamma, u) = \sum_n f_\gamma^n(\hat{y}, \hat{y}) u^n$$

$$= \sum_{\substack{n \\ \hat{y}[1] \neq \hat{y}, \ldots, \hat{y}[n-1] \neq \hat{y}, \hat{y}[n] = \hat{y}}} u^n \hat{J}_\gamma(\hat{y}, \hat{y}[1]) \hat{J}_\gamma(\hat{y}[1], \hat{y}[2]) \times \cdots$$

$$\times \hat{J}_\gamma(\hat{y}[n-2], \hat{y}[n-1]) \hat{J}_\gamma(\hat{y}[n-1], \hat{y}[n])$$

$$= \sum_{\substack{n \\ x_1, \ldots, x_n \\ \hat{y}[1] \neq \hat{y}, \ldots, \hat{y}[n-1] \neq \hat{y}}} u^n J_\gamma((0, \hat{y}), (x_1, \hat{y}[1])) J_\gamma((x_1, \hat{y}[1]), (x_2, \hat{y}[2])) \times \cdots$$

$$\times J_\gamma((x_{n-1}, \hat{y}[n-1]), (x_n, \hat{y}))$$

$$= \sum_{\substack{n \\ x_1,\ldots,x_n \\ \hat{y}[1] \neq \hat{y},\ldots,\hat{y}[n-1] \neq \hat{y}}} u^n e^{\gamma x_n} J((0, \hat{y}), (x_1, \hat{y}[1]))$$

$$\times J((x_1, \hat{y}[1]), (x_2, \hat{y}[2])) \times \cdots$$

$$\times J((x_{n-1}, \hat{y}[n-1]), (x_n, \hat{y}))$$

$$= \mathrm{E}_{(0,\hat{y})}(e^{\gamma V[T_{\hat{y}}]} u^{T_{\hat{y}}} \chi\{T_{\hat{y}} < \infty\}).$$

We introduce the random walk $\mathcal{U}$ given by $\mathcal{V}$ when $\mathcal{Z}$ returns to $\hat{y}$, where $\hat{y}$ is any fixed point in $\hat{S}$. Assume $\mathcal{U}$ starts with initial state 0. The (defective) density of the increment of this walk is given by

$$f_{\mathcal{U}}(v) = \Pr(\mathcal{V}[\mathcal{T}_{\hat{y}}] = v, \mathcal{T}_{\hat{y}} < \infty | (\mathcal{V}[0], \mathcal{Z}[0]) = (0, \hat{y})),$$

where $\mathcal{T}_{\hat{y}}$ is the first time $\mathcal{Z}$ returns to zero.

Define the transform $\Psi_{\mathcal{U}}(z) = \sum_v f_{\mathcal{U}}(v) z^v$ of $f_{\mathcal{U}}$. Let

$$\rho(x) = \Pr(\text{the first weak descending ladder point of } \mathcal{U} \text{ is at } x),$$

$$\mu(x) = \sum_{n=0}^{\infty} \Pr(\mathcal{U}_n = x \text{ and } \mathcal{U}_m > 0 \text{ for } 1 \leq m \leq n).$$

By duality (duality lemma, [5], page 395), $\mu(x)$ is equal to the probability that $x$ is a strictly ascending ladder point of $\mathcal{U}$. Let $\Psi_\rho(z) = \sum_{x \leq 0} z^x \rho(x)$ and $\Psi_\mu(z) = \sum_{x \geq 0} z^x \mu(x)$ be the transforms of $\rho$ and $\mu$, respectively.

LEMMA 3. *The radius of convergence of $\Psi_\mu(z)$ is $\sup\{u : \Psi_{\mathcal{U}}(u) \leq 1\}$.*

PROOF. We cite the following facts from [5]. $\rho + \mu = \delta_{x,0} + \mu * f_{\mathcal{U}}$. Taking transforms, we obtain

(26) $$\Psi_\mu(z) = (1 - \Psi_\rho(z))/(1 - \Psi_{\mathcal{U}}(z)).$$

Since $\mathcal{V}[n] \to \infty$ as $n \to \infty$, it follows that the return distribution is defective so $\Psi_\rho(z) < 1$ for $z > 0$. Therefore, the radius of convergence of $\Psi_\mu$ is equal to $\sup\{u : \Psi_{\mathcal{U}}(u) \leq 1\}$ [since $\Psi_{\mathcal{U}}(z)$ has radius of convergence at least one since $\mathcal{U}$ is a (defective) random variable]. □

Note that

(27) $$\begin{aligned} \Psi(\alpha + \beta, 1) &= \mathrm{E}_{(0,\hat{y})}(\exp((\alpha + \beta)V[T_{\hat{y}}])\{T_{\hat{y}} < \infty\}) \\ &= \mathrm{E}_{(0,\hat{y})}(\exp(\beta \mathcal{V}[\mathcal{T}_{\hat{y}}])\{\mathcal{T}_{\hat{y}} < \infty\}) \\ &= \sum_v f_{\mathcal{S}}(v) \exp(\beta v) = \Psi_{\mathcal{S}}(\exp(\beta)). \end{aligned}$$

LEMMA 4. *The radius of convergence of $\Psi_{\mathcal{U}}(z)$ and $\Psi_\mu(z)$ is 1.*

PROOF. $r(\hat{J}_0) = r(\hat{J}) \leq 1$ and $r(\hat{J}_\alpha) = r(\hat{\mathcal{J}}) = 1$. Since $r(\hat{J}_\gamma)$ is convex, we conclude $r(\hat{J}_{\alpha+\beta}) > 1$ for $\beta > 0$ so $R(\hat{J}_{\alpha+\beta}) < 1$ for $\beta > 0$. By the above, $\Psi_{\mathcal{S}}(\exp(\beta)) = \Psi(\alpha + \beta, 1)$. Also, by Lemma 1, for any $\beta > 0$,

$$\sup\{u : \Psi(\alpha + \beta, u) \leq 1\} = R(\hat{J}_{\alpha+\beta}) < 1 \qquad \text{so } \Psi(\alpha + \beta, 1) > 1.$$

By Lemma 3, the radius of convergence of $\Psi_{\mathcal{S}}$ is given by

$$\begin{aligned}
\sup\{u : &\Psi_{\mathcal{S}}(u) \leq 1\} \\
&= \sup\{\exp(\beta) : \Psi_{\mathcal{S}}(\exp(\beta)) \leq 1\} \\
&= \sup\{\exp(\beta) : \Psi(\alpha + \beta, 1) \leq 1\} = \exp(0) = 1.
\end{aligned}$$

Hence, the radius of convergence of $\Psi_{\mathcal{S}}$ is one, and by Lemma 3, the radius of convergence of $\Psi_\mu$ is one as well. □

Consider the Markov chain $A(\ell)$ having kernel $\mathcal{A}$, which represents the age of the ladder height process at time $\ell$. The probability $\mu(\ell)$ is the probability that the age at time $\ell$ is zero; that is, $\mu(\ell) = P(A(\ell) = 0) = \mathcal{A}_{0,0}^\ell$. Next $(\mathcal{A}_{0,0}^\ell)^{1/\ell} \to r(\mathcal{A})$, the spectral radius of $\mathcal{A}$, so we conclude $\mu(\ell)^{1/\ell} \to 1$. Let $m(\hat{y})$ be the expected number of visits to $(0, \hat{y})$ by $(\mathcal{V}, \mathcal{Z})$ starting from the point $(0, \hat{y})$. Clearly,

$$m(\hat{y}) \geq \mathcal{G}((0, \hat{y}), (\ell, \hat{y})) \geq \mu(\ell) m(\hat{y}).$$

Since $\mu(\ell)^{1/\ell} \to 1$, this means $\lim_{\ell \to \infty} (\mathcal{G}((0, \hat{y}), (\ell, \hat{y})))^{1/\ell} = 1$.

5.2. *Proof of* (20).

LEMMA 5. *Under assumptions* A0–A6,

$$\mathcal{G}((0, \hat{z}), (\ell + 1, \hat{y})) / \mathcal{G}((0, \hat{z}), (\ell, \hat{y})) \to 1.$$

For modified Jackson networks the asymptotics of $\mathcal{G}((0,0), (\ell, 0))$ are calculated explicitly in Proposition 1 and can be directly shown to satisfy the conclusion of Lemma 5.

PROOF. Let $\mathcal{X}_k = \mathcal{V}[k] - \mathcal{V}[k-1]$. By A1,

$$P(\mathcal{V}[n] \geq \ell | \mathcal{Z}[0] = \hat{z}, \mathcal{V}[0] = 0)$$
$$\leq e^{-t\ell} \mathrm{E}\left(\exp\left(t \sum_{k=1}^n \mathcal{X}_k\right) \Big| \mathcal{Z}[0] = \hat{z}, \mathcal{V}[0] = 0\right)$$

$$\leq e^{-t\ell} \mathrm{E}\left(\mathrm{E}\left(\exp\left(t\sum_{k=1}^{n}\mathcal{X}_k\right)\bigg|\mathcal{Z}[k], k=1,2,\ldots,n:\mathcal{Z}[0]=\hat{z},\mathcal{V}[0]=0\right)\bigg|\right.$$
$$\left.\mathcal{Z}[0]=\hat{z},\mathcal{V}[0]=0\right)$$
$$\leq e^{-t\ell}M^n.$$

Next, for any $\kappa > 0$,

$$\sum_{n=0}^{\kappa\ell} P(\mathcal{Z}[n]=\hat{y},\mathcal{V}[n]=\ell+1|\mathcal{V}[0]=0,\mathcal{Z}[0]=\hat{z})$$
$$\leq \sum_{n=0}^{[\kappa\ell]} P(\mathcal{V}[n]\geq\ell)$$
$$\leq \sum_{n=0}^{[\kappa\ell]} e^{-t\ell}M^n \leq e^{-t\ell}M^{\kappa\ell}/(M-1),$$

where $[x]$ is the integer part of $x$. We pick $\kappa$ sufficiently small so that $s = \exp(-t)M^\kappa < 1$. Hence, we can decompose

$$(28) \quad \mathcal{G}((0,\hat{z}),(\ell+1,\hat{y})) = \sum_{n=0}^{\infty} P(\mathcal{Z}[n]=\hat{y},\mathcal{V}[n]=\ell+1|\mathcal{Z}[0]=\hat{z},\mathcal{V}[0]=0)$$

into a main part

$$(29) \quad \sum_{n=\kappa\ell}^{\infty} P(\mathcal{Z}[n]=\hat{y},\mathcal{V}[n]=\ell+1|\mathcal{Z}[0]=\hat{z},\mathcal{V}[0]=0)$$

and a negligible part. We take negligible to mean that a term like $\exp(-s\ell)$ is negligible compared to $\mathcal{G}((0,\hat{y}),(\ell,\hat{y}))$ as $\ell\to\infty$. Such a term is negligible because $\lim_{\ell\to\infty}(\mathcal{G}((0,\hat{y}),(\ell,\hat{y})))^{1/\ell}=1$; that is, the potential does not die out exponentially fast. By assumption A2, $\mathcal{J}^m((0,\hat{z}),(0,\hat{y})) > 0$ for some $m$. It follows that

$$(30) \quad \mathcal{G}((0,\hat{z}),(\ell,\hat{y})) \geq \mathcal{J}^m((0,\hat{z}),(0,\hat{y}))\mathcal{G}((0,\hat{y}),(\ell,\hat{y}))$$

modulo the negligible probability that $\mathcal{V}$ exceeds $\ell$ in $m$ steps. Consequently, if a term is negligible compared to $\mathcal{G}((0,\hat{y}),(\ell,\hat{y}))$, it is negligible compared to $\mathcal{G}((0,\hat{z}),(\ell,\hat{y}))$.

Without loss of generality, we can assume there is a minimum probability the Markov additive process $(\mathcal{V}[n],\mathcal{Z}[n])$ stays put during any transition. By assumption A6, there is a fixed minimum probability that the Markov additive process $(\mathcal{V}[n],\mathcal{Z}[n])$ makes a transition from $(0,\hat{z})$ to $(1,\hat{z})$ after $N$ steps. Hence,

$$\min\{\mathcal{J}^N((0,\hat{z}),(0,\hat{z})),\mathcal{J}^N((0,\hat{z}),(1,\hat{z}))\} \geq \delta \qquad \text{uniformly in } \hat{z}.$$

We can use the Bernoulli part decomposition developed in [3] and [13] to represent

$$\mathcal{V}[n] = \mathcal{U}[n] + \sum_{k=1}^{N_n} L_k,$$

where $L_1, L_2, \ldots$ are i.i.d. Bernoulli random variables independent of $\{(\mathcal{U}[n], N_n); n = 1, 2, \ldots\}$ such that $P(L_1 = 0) = P(L_1 = 1) = 1/2$ and where $N_n$ is a Bernoulli random variable with mean $nb$ and variance $nb(1-b)$, where $b = \delta/N$. $N_n$ and $U[n]$ are dependent. In effect, we have a probability $\delta$ of picking up a Bernoulli step every $N$ transitions.

We can represent (29) as

$$\sum_{n=\kappa\ell}^{\infty} P\bigg(\mathcal{Z}[n] = \hat{y}, \mathcal{U}[n] + \sum_{k=1}^{N_n} L_k = \ell + 1\bigg)$$

(31)
$$= \sum_{n=\kappa\ell}^{\infty} \sum_{m=1}^{n} \sum_{x=1}^{m} P(\mathcal{Z}[n] = \hat{y}, U[n] = \ell - x, N_n = m)$$

$$\times P\bigg(\sum_{k=1}^{m} L_k = x + 1\bigg).$$

Pick $\varepsilon$ to be small and such that $0 < \varepsilon < b$. For any $n$, we can bound the large deviation probability

$$P(|N_n - nb| > \varepsilon n) \leq \exp(-\Lambda_1 n) \qquad \text{where } \Lambda_1 > 0,$$

as in [22], Example 1.15. Therefore, we can decompose (31) into a main part

(32)
$$\sum_{n=\kappa\ell}^{\infty} \sum_{m:\,|n-bn|\leq \varepsilon n} \sum_{x=1}^{m} P(\mathcal{Z}[n] = \hat{y}, U[n] = \ell - x, N_n = m)$$

$$\times P\bigg(\sum_{k=1}^{m} L_k = x + 1\bigg)$$

and a negligible part

$$\sum_{n=\kappa\ell}^{\infty} \sum_{m:\,|n-bn|> \varepsilon n} \sum_{x=1}^{m} P(\mathcal{Z}[n] = \hat{y}, U[n] = \ell - x, N_n = m) P\bigg(\sum_{k=1}^{m} L_k = x + 1\bigg)$$

$$\leq \sum_{n=\kappa\ell}^{\infty} \sum_{m:\,|n-bn|> \varepsilon n} P(N_n = m)$$

$$\leq \sum_{n=\kappa\ell}^{\infty} \exp(-\Lambda_1 n)$$

$$= \exp(-\kappa\ell\Lambda_1)/(1 - \exp(-\Lambda_1)).$$

We can now expand (32) into a main part

$$\sum_{n=\kappa\ell}^{\infty} \sum_{m:\,|m-bn|\leq\varepsilon n} \sum_{x:\,|x-m/2|\leq\varepsilon m} P(\mathcal{Z}[n]=\hat{y}, U[n]=\ell-x, N_n=m)$$
(33)
$$\times P\left(\sum_{k=1}^{m} L_k = x+1\right)$$

and a negligible part

$$\sum_{n=\kappa\ell}^{\infty} \sum_{m:\,|m-bn|\leq\varepsilon n} \sum_{x:\,|x-m/2|>\varepsilon m} P(\mathcal{Z}[n]=\hat{y}, U[n]=\ell-x, N_n=m)P$$
$$\times \left(\sum_{k=1}^{m} L_k = x+1\right)$$
$$\leq \sum_{n=\kappa\ell}^{\infty} \sum_{m:\,|m-bn|\leq\varepsilon n} P(\mathcal{Z}[n]=\hat{y}, N_n=m)\exp(-\Lambda_2 m)$$
$$\leq \sum_{n=\kappa\ell}^{\infty} \sum_{m\geq(b-\varepsilon)n} \exp(-\Lambda_2 m)$$
$$\leq \sum_{n=\kappa\ell}^{\infty} \exp(-\Lambda_2(b-\varepsilon)n)/(1-\exp(-\Lambda_2))$$
$$\leq \exp(-\Lambda_2(b-\varepsilon)\kappa\ell)/((1-\exp(-\Lambda_2(b-\varepsilon)))(1-\exp(-\Lambda_2))),$$

where the constant $\Lambda_2 > 0$ is given in [22], Example 1.15.

We can now focus on (33),

$$P\left(\sum_{k=1}^{m} L_k = x+1\right) = \binom{m}{x+1}\frac{1}{2^m} = P\left(\sum_{k=1}^{m} L_k = x\right)\frac{m-x}{x+1}.$$

For $x$ such that $|x-m/2|\leq\varepsilon m$, we have

$$\frac{1/2-\varepsilon}{1/2+\varepsilon+1/m} \leq \frac{m-x}{x+1} \leq \frac{1/2+\varepsilon}{1/2-\varepsilon}.$$

If, in addition, $|m-bn|\leq\varepsilon n$, we can simplify the left-hand side of the above:

$$\frac{1/2-\varepsilon}{1/2+\varepsilon+1/(n(b+\varepsilon))} \leq \frac{m-x}{x+1}.$$

Finally, when $n\geq\kappa\ell$, we get

$$f_1 \equiv \frac{1/2-\varepsilon}{1/2+\varepsilon+1/((b+\varepsilon)\kappa\ell)} \leq \frac{m-x}{x+1} \leq \frac{1/2+\varepsilon}{1/2-\varepsilon} \equiv f_2.$$

We can therefore bound (33) above and below by multiplying $f_2$ and $f_1$, respectively, by

$$\sum_{n=\kappa\ell}^{\infty} \sum_{m\,:\,|m-bn|\leq\varepsilon n} \sum_{x\,:\,|x-m/2|\leq\varepsilon m} P\bigg(\mathcal{Z}[n]=\hat{y},U[n]=\ell-x,N_n=m\bigg)$$
(34)
$$\times P\bigg(\sum_{k=1}^{m} L_k = x\bigg).$$

If we now add the negligible terms back into (34), we get that (31) is bounded (modulo negligible terms) above and below by multiplying $f_2$ and $f_1$, respectively, by

$$\sum_{n=0}^{\infty} P\bigg(\mathcal{Z}[n]=\hat{y},\mathcal{U}[n]+\sum_{k=1}^{N_n}L_k=\ell\bigg) = \mathcal{G}((0,\hat{z}),(\ell,\hat{y})).$$

Hence,

$$\frac{1/2-\varepsilon}{1/2+\varepsilon} \leq \liminf_{\ell\to\infty}\frac{\mathcal{G}((0,\hat{z}),(\ell+1,\hat{y}))}{\mathcal{G}((0,\hat{z}),(\ell,\hat{y}))} \leq \limsup_{\ell\to\infty}\frac{\mathcal{G}((0,\hat{z}),(\ell+1,\hat{y}))}{\mathcal{G}((0,\hat{z}),(\ell,\hat{y}))} \leq \frac{1/2+\varepsilon}{1/2-\varepsilon}$$

and since $\varepsilon$ can be taken arbitrarily small, the result follows. $\square$

5.3. *Proof of* (21)—*uniform integrability.*

PROPOSITION 4. *Under assumptions* A0–A7

$$\mathcal{G}((t,\hat{z}),(\ell,\hat{\sigma})) \sim \mathcal{G}((0,\hat{\sigma}),(\ell,\hat{\sigma})) \qquad as\ \ell\to\infty.$$

PROOF. Define

$$B_{\ell_n}(s,\hat{x}) := \mathcal{G}((s,\hat{x}),(\ell_n,\hat{\sigma}))/\mathcal{G}((0,\hat{\sigma}),(\ell_n,\hat{\sigma})).$$

Take a subsequence $\ell_n$ such that $B_{\ell_n}(s,\hat{x})$ converges to $B(s,\hat{x})$ for all $(s,\hat{x})$ as $\ell_n\to\infty$. By Lemma 5, $B(s,\hat{x})=B(\hat{x})$.

We can write

(35) $$B_{\ell_n}(t,\hat{z}) = \sum_{(s,\hat{x})} \mathcal{J}((t,\hat{z}),(s,\hat{x}))B_{\ell_n}(s,\hat{x}) + \frac{\mathcal{J}((t,\hat{z}),(\ell_n,\hat{\sigma}))}{\mathcal{G}((0,\hat{\sigma}),(\ell_n,\hat{\sigma}))}.$$

By assumption A.1, $\mathcal{J}((t,\hat{z}),(\ell,\hat{\sigma}))$ above decays exponentially fast as $\ell\to\infty$ so the second term on the right-hand side of (35) tends to zero as $\ell\to\infty$. Moreover,

$$\mathcal{J}((t,\hat{z}),(s,\hat{x})) = \hat{\mathcal{J}}(\hat{z},\hat{x})\tilde{\mathcal{J}}_{(\hat{z},\hat{x})}(s-t),$$

where $\tilde{\mathcal{J}}_{(\hat{z},\hat{x})}(s-t)$ is the probability the additive transition equals $s-t$ given there has been a transition from $\hat{z}$ to $\hat{x}$,

$$B_{\ell_n}(s,\hat{x}) = \frac{\mathcal{G}((s,\hat{x}),(\ell_n,\hat{\sigma}))}{\mathcal{G}((0,\hat{\sigma}),(\ell_n,\hat{\sigma}))} = \frac{\mathcal{G}((s,\hat{x}),(\ell_n,\hat{\sigma}))}{\mathcal{G}((s,\hat{\sigma}),(\ell_n,\hat{\sigma}))} \frac{\mathcal{G}((s,\hat{\sigma}),(\ell_n,\hat{\sigma}))}{\mathcal{G}((0,\hat{\sigma}),(\ell_n,\hat{\sigma}))}$$

and the first fraction is bounded by $f(\hat{x})$ by A7. Moreover,

$$\sum_{\hat{x}} \hat{\mathcal{J}}(\hat{z},\hat{x}) f(\hat{x}) < \infty$$

for all $\hat{z}$.

If $p((0,\hat{\sigma}),(t,\hat{\sigma}))$ is the probability of ever going from $(0,\hat{\sigma})$ to $(t,\hat{\sigma})$, then

$$p((0,\hat{\sigma}),(\ell_n,\hat{\sigma})) \geq p((0,\hat{\sigma}),(s,\hat{\sigma})) p((s,\hat{\sigma}),(\ell_n,\hat{\sigma})).$$

Hence,

$$\frac{\mathcal{G}((s,\hat{\sigma}),(\ell_n,\hat{\sigma}))}{\mathcal{G}((0,\hat{\sigma}),(\ell_n,\hat{\sigma}))} = \frac{p((s,\hat{\sigma}),(\ell_n,\hat{\sigma}))}{p((0,\hat{\sigma}),(\ell_n,\hat{\sigma}))} \leq \frac{1}{p((0,\hat{\sigma}),(s,\hat{\sigma}))}.$$

$1/p((0,\hat{\sigma}),(s,\hat{\sigma}))$ increases at a subexponential rate in $s$ and, hence, is integrable in $s$ with respect to $\tilde{\mathcal{J}}_{(\hat{z},\hat{x})}(s-t)$ since $\tilde{\mathcal{J}}_{(\hat{z},\hat{x})}(s-t)$ decays exponentially fast in $s$ uniformly in $\hat{z},\hat{x}$ by A1.

We conclude that $B_{\ell_n}(s,\hat{x})$ is bounded uniformly in $\ell_n$ by $f(\hat{x})/p((0,\hat{\sigma}),(s,\hat{\sigma}))$ and

$$\sum_{(s,\hat{x})} \mathcal{J}((t,\hat{z}),(s,\hat{x})) f(\hat{x})/p((0,\hat{\sigma}),(s,\hat{\sigma})) < \infty.$$

Since $B_{\ell_n}(s,\hat{x})$ converges pointwise to $B(\hat{x})$, dominated convergence implies that $B_{\ell_n}(s,\hat{x})$ converges in $L^1$ relative to the measure $\mathcal{J}((t,\hat{z}),(s,\hat{x}))$ and that

$$B(t,\hat{z}) = \sum_{(s,\hat{x})} \mathcal{J}((t,\hat{z}),(s,\hat{x})) B(s,\hat{x})$$

or

$$B(\hat{z}) = \sum_{\hat{x}} \hat{\mathcal{J}}(\hat{z},\hat{x}) B(\hat{x})$$

by Lemma 5. This means $B(\hat{z})$ is harmonic for $\hat{\mathcal{J}}$. By hypothesis, the constant functions are the only harmonic functions for $\hat{\mathcal{J}}$ and $B(\hat{\sigma}) = 1$ so $B(s,\hat{\sigma}) = 1$ for all $(s,\hat{z})$. Since this limit is independent of the subsequence $\ell_n$, it follows that

$$\lim_{\ell \to \infty} \mathcal{G}((t,\hat{z}),(\ell,\hat{\sigma}))/\mathcal{G}((0,\hat{\sigma}),(\ell,\hat{\sigma})) = 1. \qquad \square$$

5.4. *Proof of* (22).

PROPOSITION 5.

(36) $$\frac{\mathcal{G}(z,(\ell,\hat{y}))}{\mathcal{G}((0,\hat{y}),(\ell,\hat{y}))} \leq \frac{1}{p((0,\hat{y}),(0,\hat{\sigma}))^2 p((0,\hat{\sigma}),(0,\hat{y}))} \frac{\mathcal{G}(z,(\ell,\hat{\sigma}))}{\mathcal{G}((0,\hat{\sigma}),(\ell,\hat{\sigma}))}.$$

PROOF. The number of visits to $(\ell,\hat{\sigma})$ is greater than the number of visits to $(\ell,\hat{y})$ followed by a visit to $(\ell,\hat{\sigma})$,

$$\mathcal{G}(z,(\ell,\hat{\sigma})) \geq \mathcal{G}(z,(\ell,\hat{y}))p((0,\hat{y}),(0,\hat{\sigma})),$$

and the number of visits to $(\ell,\hat{y})$ is greater than the number of visits to $(\ell,\hat{y})$ that must first go through $(0,\hat{\sigma})$ and then must go through $(0,\hat{\sigma})$ before going to $(\ell,\hat{y})$,

$$\mathcal{G}((0,\hat{y}),(\ell,\hat{y})) \geq p((0,\hat{y}),(0,\hat{\sigma}))\mathcal{G}((0,\hat{\sigma}),(\ell,\hat{\sigma}))p((0,\hat{\sigma}),(0,\hat{y})).$$

The inequality follows by dividing the above two inequalities □

PROPOSITION 6. *Under assumptions* A0–A7,

$$\mathcal{G}((s,\hat{z}),(\ell,\hat{y})) \sim \mathcal{G}((0,\hat{y}),(\ell,\hat{y})) \qquad as\ \ell \to \infty.$$

PROOF. Define

$$C_{\ell_n}(s,\hat{x}) := \mathcal{G}((s,\hat{x}),(\ell_n,\hat{y}))/\mathcal{G}((0,\hat{y}),(\ell_n,\hat{y})).$$

Take a subsequence $\ell_n$ such that $C_{\ell_n}(s,\hat{x})$ converges to $C(s,\hat{x})$ for all $(s,\hat{x})$ as $\ell_n \to \infty$. By Lemma 5, $C(s,\hat{x}) = C(\hat{x})$. By Proposition 5,

$$C_{\ell_n}(s,\hat{x}) \leq \frac{1}{p((0,\hat{y}),(0,\hat{\sigma}))^2 p((0,\hat{\sigma}),(0,\hat{y}))} \frac{\mathcal{G}((s,\hat{x}),(\ell_n,\hat{\sigma}))}{\mathcal{G}((0,\hat{\sigma}),(\ell_n,\hat{\sigma}))}.$$

The sequence of functions in $(s,\hat{x})$ indexed by $\ell_n$ on the right-hand side above is uniformly integrable with respect to the measure $\mathcal{J}((t,\hat{z}),(s,\hat{x}))$, so $C_{\ell_n}(s,\hat{x})$ is uniformly integrable as well.

Since

$$C_{\ell_n}(t,\hat{z}) = \sum_{(s,\hat{x})} \mathcal{J}((t,\hat{z}),(s,\hat{x}))C_{\ell_n}(s,\hat{x}) + \frac{\mathcal{J}((t,\hat{z}),(\ell_n,\hat{y}))}{\mathcal{G}((0,\hat{y}),(\ell_n,\hat{y}))},$$

we can take the limit as $\ell_n \to \infty$ to get

$$C(\hat{x}) = \sum_{\hat{x}} \hat{\mathcal{J}}(\hat{z},\hat{x})C(\hat{x}).$$

Again, we conclude $C(\hat{x}) = 1$ and this proves the result. □

Define the time reversed kernel by

$$\mathcal{J}^*(w, z) = \varphi(\hat{z})\mathcal{J}(z, w)/\varphi(\hat{w}).$$

$\mathcal{J}^*$ is the kernel of the time reversed process $\mathcal{W}^* = (\mathcal{V}^*, \mathcal{Z}^*)$. Define the Green's function, $\mathcal{G}^*$ by

$$\mathcal{G}^*(z, w) \equiv \mathcal{G}^*((z_1, \hat{z}), (w_1, \hat{w})) = \sum_{k=0}^{\infty} (\mathcal{J}^*)^n(z, w),$$

$\mathcal{Z}^*$ has kernel $\hat{\mathcal{J}}^*(\hat{x}, \hat{z}) = \varphi(\hat{z})\hat{\mathcal{J}}(\hat{z}, \hat{x})/\varphi(\hat{x})$ and $\hat{\mathcal{J}}^*$ has unique (up to constants) invariant measure $\varphi$ and harmonic function 1.

PROPOSITION 7. *Under assumptions* A0–A6 *and* A7*,

$$\mathcal{G}((0, \hat{y}), (\ell, \hat{x}))/\mathcal{G}((0, \hat{y}), (\ell, \hat{y})) \to \frac{\varphi(\hat{x})}{\varphi(\hat{y})} \qquad as\ \ell \to \infty.$$

PROOF. It is easy to see $-\mathcal{W}^*[n]$ satisfies assumptions A0–A6. Moreover, A7 holds for $-\mathcal{W}^*[n]$ because A7* holds for $\mathcal{W}$. Consequently,

$$\mathcal{G}^*((0, \hat{x}), (-\ell, \hat{y}))/\mathcal{G}^*((0, \hat{y}), (-\ell, \hat{y})) \to 1 \qquad as\ \ell \to \infty$$

by Proposition 4 since the constant functions are the unique harmonic functions for $\hat{\mathcal{J}}^*$. By time reversal, this means

$$\left(\frac{\varphi(\hat{y})}{\varphi(\hat{x})}\mathcal{G}((0, \hat{y}), (\ell, \hat{x}))\right) \Big/ \left(\frac{\varphi(\hat{y})}{\varphi(\hat{y})}\mathcal{G}((0, \hat{y}), (\ell, \hat{y}))\right) \to 1 \qquad as\ \ell \to \infty.$$

This gives the result. □

Putting Propositions 4 and 7 together we get the following:

THEOREM 2. *Under assumptions* A0–A7 *and* A7*

$$\mathcal{G}((s, \hat{z}), (\ell, \hat{x}))/\mathcal{G}((u, \hat{w}), (\ell, \hat{y})) \to \frac{\varphi(\hat{x})}{\varphi(\hat{y})} \qquad as\ \ell \to \infty$$

*for fixed* $s, u, \hat{w}, \hat{x}, \hat{y}$ *and* $\hat{z}$.

Using Theorem 2, plus the change of measure from $J$ to $\mathcal{J}$, we get the following:

THEOREM 3. *Under assumptions* A0–A7 *and* A7*,

$$G((s, \hat{z}), (\ell, \hat{x}))/G((u, \hat{w}), (\ell, \hat{y})) \to \frac{\hat{h}(\hat{x})\varphi(\hat{x})}{\hat{h}(\hat{y})\varphi(\hat{y})} \frac{\hat{h}(\hat{w})}{\hat{h}(\hat{z})} \qquad as\ \ell \to \infty$$

*for fixed* $s, u, \hat{w}, \hat{x}, \hat{y}$ *and* $\hat{z}$.

5.5. *Assumptions* A0–A6 *are insufficient for* (21). This section gives an example satisfying all of the assumptions except for the uniform integrability assumptions A7 and A7*, but (21) does not hold. From the discussion near (25), it suffices to give a transition kernel $\hat{\mathcal{J}}$ such that (24) does not hold. Consider a chain on $\hat{S} = \{\ldots, -2, -1, 0, 1, 2, \ldots\}$ with kernel $\hat{\mathcal{J}}$ defined as follows:

$$\hat{\mathcal{J}}(n, n+1) = \hat{\mathcal{J}}(n, n-1) = 1/2 \quad \text{for } n \geq 1$$
$$\hat{\mathcal{J}}(n, n+1) = 1/3 \quad \text{and} \quad \hat{\mathcal{J}}(n, n-1) = 2/3 \quad \text{for } n \leq -1$$
$$\hat{\mathcal{J}}(0, n) = f(n) \quad \text{for } n \geq -1 \text{ and } f(-1) > 0,$$

where $f$ is the probability mass function such that $\sum_{n=0}^{\infty} nf(n) = \infty$. We could modify the kernel so $\hat{\mathcal{J}}(n,n) > \varepsilon$ for all $n$ and this will not alter the following conclusions.

This chain is clearly irreducible and transient to minus infinity. Nevertheless, it has spectral radius one. Starting from 0, the chain can go to 1 and then return to 0 after $2\ell - 1$ steps without hitting 0 again. We saw this has probability of order $1/(2\ell)^{3/2}$ in the preceding section. Hence,

$$\lim_{\ell \to \infty} (\hat{\mathcal{J}}^{2\ell}(0,0))^{1/2\ell} = 1.$$

This chain only has constant positive harmonic functions because for $n > 1$, a harmonic function $h$ must satisfy $h(n) = \frac{1}{2}h(n+1) + \frac{1}{2}h(n-1)$. This means the harmonic function must be positive and linear on $[0, \infty)$. On the other hand, we must have $h(0) = \sum_{n=0}^{\infty} h(n)f(n)$ and this is infinity unless $h$ is constant by the construction of $f$. Hence, no harmonic functions other than constants can exist.

In spite of all this, Proposition 4 fails with $\hat{y} = 0$. Pick a starting point $\hat{x} = -n$ with $n$ even. Once the chain hits 0 the first time after $2k$ steps, the probability of hitting 0 again in $2\ell - 2k$ steps is asymptotically the same as returning to 0 in $2\ell$ steps since this later probability is of order $1/(2\ell)^{3/2}$. It therefore follows that

$$\frac{\hat{\mathcal{J}}^{2\ell}(-n, 0)}{\hat{\mathcal{J}}^{2\ell}(0, 0)} \to P(\text{hit 0 starting from } -n).$$

We bound the probability of hitting 0 from $-n$ by considering 0 to be absorbing and we see the return distribution falls off geometrically fast in $n$. Hence, the probability of a return to 0 can be picked arbitrarily small by picking $n$ large. Hence, $P(\text{hit 0 starting from } -n)$ is strictly less than one. Hence, the ratio limit theorem fails. This is not unexpected since jumps are unbounded and A7 fails.

**6. Asymptotics of queueing networks.** As in [8, 14], we consider a Markov additive process $W^\infty$ with transition kernel $K^\infty$ on $S^\infty$ and a boundary ▲ with an edge denoted by $\triangle = S \cap$ ▲ such that the probability transition kernel $K$ of the chain $W$ agrees with $K^\infty$ within $S^\infty \setminus$ ▲; that is, $K(x,C) = K^\infty(x,C)$ if $x \in S$ and $C \subset S^\infty \setminus$ ▲.

We assume $W$ is stable and we are interested in the asymptotics of the steady state probability $\pi$ of $W$ as the additive component gets large. We have assumed the existence of a positive harmonic function $h(w) \equiv h(w) = \exp(\alpha w_1)\hat{h}(\hat{w})$ for $K^\infty$ on $S^\infty$. We have constructed the $h$-transform $\mathcal{K}^\infty(z,w) := K^\infty(z,w) \times h(w)/h(z)$. Let $\hat{\mathcal{K}}^\infty$ denote the Markovian part of the transition kernel $\mathcal{K}^\infty$. We assume $(\mathcal{W}_1^\infty, \hat{\mathcal{W}}^\infty) \equiv (\mathcal{V}, \mathcal{Z})$ satisfies A0–A7, A7* and B1–B4 given below:

B1a. Define Green's function,

$$\mathcal{G}(z,y) \equiv \mathcal{G}((z_1,\hat{z}),(y_1,\hat{y})) = \sum_{k=0}^\infty (\mathcal{K}^\infty)^n(z,y).$$

We assume that for some fixed state $\hat{\sigma}$,

$\dfrac{\mathcal{G}((0,\hat{w}),(\ell,\hat{\sigma}))}{\mathcal{G}((0,\hat{\sigma}),(\ell,\hat{\sigma}))}$ is bounded uniformly in $\hat{w}$ for $\ell$ sufficiently large.

B1b. We also assume

$\dfrac{\mathcal{G}(w,(\ell,\hat{\sigma}))}{\mathcal{G}((0,\hat{\sigma}),(\ell,\hat{\sigma}))}$ is bounded uniformly in $w \in$ ▲ for $\ell$ sufficiently large.

In fact, the uniform boundedness need only be checked on that subset of ▲ which can be reached in a transition of $W^\infty$ from $S^\infty \setminus$ ▲ into ▲.

B2. There exists a subset $C \subseteq \triangle$ such that $\pi(C) > 0$ and such that, for $z \in C$, $P_z(\mathcal{T}_\blacktriangle^\infty = \infty) > 0$, where $\mathcal{T}_\blacktriangle^\infty$ is the first return time to ▲ by $\mathcal{W}^\infty$.

B3. For each $\hat{y}$, there is an associated integer $L(\hat{y})$ such that $(\ell,\hat{y}) \cap \triangle = \varnothing$ if $\ell \geq L(\hat{y})$.

B4. $\lambda(x) \equiv \pi(x)h(x)\chi\{x \in \triangle\}$ is a finite measure.

LEMMA 6. *If assumption* B1 *holds, then for any state $\hat{y}$ and for $\ell$ sufficiently large,*

$\dfrac{\mathcal{G}(z,(\ell,\hat{y}))}{\mathcal{G}((0,\hat{y}),(\ell,\hat{y}))}$ *is bounded uniformly in $z \in$ ▲ for $\ell$ sufficiently large*

*and*

$\dfrac{\mathcal{G}((0,\hat{z}),(\ell,\hat{y}))}{\mathcal{G}((0,\hat{y}),(\ell,\hat{y}))}$ *is bounded uniformly in $\hat{z}$ for $\ell$ sufficiently large.*

PROOF. By Proposition 5,
$$\frac{\mathcal{G}(z,(\ell,\hat{y}))}{\mathcal{G}((0,\hat{y}),(\ell,\hat{y}))} \leq C \frac{\mathcal{G}(z,(\ell,\hat{\sigma}))}{\mathcal{G}((0,\hat{\sigma}),(\ell,\hat{\sigma}))} \qquad \text{where } C \text{ is a constant.}$$

The right-hand side is uniformly bounded in $z \in \blacktriangle$ for $\ell$ sufficiently large by condition B1b.

The same argument works to prove the second inequality, except condition B1a is used. $\square$

THEOREM 4. *Under assumptions* A0–A7, B1–B4, *for any states $\hat{y}$ and $\hat{\sigma}$ such that $\varphi(\hat{\sigma}) > 0$,*
$$e^{\alpha \ell}\pi(\ell,\hat{y}) \sim f\hat{h}^{-1}(\hat{y})\varphi(\hat{y})\frac{\mathcal{G}((0,\hat{\sigma}),(\ell,\hat{\sigma}))}{\varphi(\hat{\sigma})},$$
*where $f = \sum_{z\in\triangle}\pi(z)h(z)P_z(\mathcal{T}_{\blacktriangle}^\infty = \infty)$.*

Note $f$ is positive by assumption B2 and can be obtained by fast simulation. Also note that if $\hat{\mathcal{K}}^\infty$ is positive recurrent, the *steady state theorem* in [8, 14] implies
$$\frac{\mathcal{G}((0,\hat{\sigma}),(\ell,\hat{\sigma}))}{\varphi(\hat{\sigma})} \to \frac{1}{\tilde{\mu}} \qquad \text{where } \tilde{\mu} \text{ is the mean drift of } \tilde{\mathcal{W}}_1^\infty.$$

PROOF OF THEOREM 4. As in the steady state theorem in [8], for $\ell > L(\hat{y})$, the steady state probability of $(\ell,\hat{y})$ is given by
$$\pi(\ell,\hat{y}) = \sum_{z\in\triangle}\pi(z)\mathrm{E}_z\left(\sum_{n=0}^{T_{\blacktriangle}^\infty-1}\chi\{W^\infty[n]=(\ell,\hat{y})\}\right)$$
$$= \sum_{z\in\triangle}\pi(z)\mathrm{E}_z(N_Y^{\blacktriangle}(\ell)),$$
where $N_Y^{\blacktriangle}(\ell)$ denotes the number of visits by $W^\infty$ to $(\ell,\hat{y})$ before $T_{\blacktriangle}^\infty$; that is, before $W^\infty$ returns to $\blacktriangle$. Also let $N(\ell,\hat{y})$ denote the number of visits by $W^\infty$ to $(\ell,\hat{y})$.

By the change of measure induced by the twist, we get
$$(37) \qquad \alpha^\ell \pi(\ell,\hat{y}) = \hat{h}^{-1}(\hat{y})\sum_{z\in\triangle}\pi(z)h(z)\mathrm{E}_z(\mathcal{N}^{\blacktriangle}(\ell,\hat{y})),$$
where $\mathcal{N}^{\blacktriangle}(\ell,\hat{y})$ is the number of visits by the twisted process. If $\mathcal{N}(\ell,\hat{y})$ is the number of visits to $(\ell,\hat{y})$ by $\mathcal{W}^\infty$, then $\mathcal{N}^{\blacktriangle}(\ell,\hat{y})$ agrees with $\mathcal{N}(\ell,\hat{y})$ if $\mathcal{W}^\infty$ never hits $\blacktriangle$.

Let $m(\ell,\hat{y}) = \mathcal{G}((0,\hat{y}),(\ell,\hat{y}))$. We now investigate
$$(38) \qquad \sum_{z\in\triangle}\pi(z)h(z)\frac{1}{m(\ell,\hat{y})}\mathrm{E}_z(\mathcal{N}^{\blacktriangle}(\ell,\hat{y}))$$

as $\ell \to \infty$.

Let $\mathcal{T}_b^\infty = \min\{n \geq 0 : \tilde{\mathcal{W}}_1^\infty[n] \geq b\}$,

$$\sum_{z \in \triangle} \lambda(z) P_z(\mathcal{T}_b^\infty < \mathcal{T}_\blacktriangle^\infty < \infty)$$

$$= \sum_{z \in \triangle} \pi(z) \sum_{w \in \triangle} P_z(T_b^\infty < T_\blacktriangle^\infty < \infty, W(T_\blacktriangle^\infty) = w) h(w)$$

$$= \sum_{w \in \triangle} \pi(w) h(w) \sum_{z \in \triangle} P_w(T_b^* < T_\blacktriangle^* < \infty, W(T_\blacktriangle^*) = z)$$

where "$*$" indicates time reversal with respect to $\pi$

$$\leq \sum_{w \in \triangle} \pi(w) h(w) P_w(T_b^* < T_\blacktriangle^* < \infty)$$

$$\to 0$$

since $\lambda(z) \equiv \pi(z) h(z) \chi\{z \in \triangle\}$ is a finite measure by hypothesis and $\pi(W_1^* \geq b) \to 0$ as $b \to \infty$. Pick $b$ sufficiently big so that

$$\sum_{z \in \triangle} \lambda(z) P_z(\mathcal{T}_b^\infty < \mathcal{T}_\blacktriangle^\infty < \infty) < \varepsilon.$$

Moreover, by Lemma 6, we can also pick $b$ sufficiently large that

$$\sum_{z \in \triangle, \tilde{z}_1 \geq b} \pi(z) h(z) \frac{1}{m(\ell, \hat{y})} \mathrm{E}_z(\mathcal{N}^\blacktriangle(\ell, \hat{y}))$$

(39)

$$\leq \sum_{z \in \triangle, \tilde{z}_1 \geq b} \pi(z) h(z) \frac{1}{m(\ell, \hat{y})} \mathrm{E}_z(\mathcal{N}(\ell, \hat{y})) < \varepsilon$$

for $\ell$ sufficiently large.

For $\ell > b$, we first remark that

(40)
$$\sum_{z \in \triangle, \tilde{z}_1 \leq b} \pi(z) h(z) \frac{1}{m(\ell, \hat{y})} \mathrm{E}_z(\mathcal{N}^\blacktriangle(\ell, \hat{y}))$$

$$= \sum_{z \in \triangle, \tilde{z}_1 \leq b} \pi(z) h(z) \frac{1}{m(\ell, \hat{y})} \mathrm{E}_z(\mathcal{N}^\blacktriangle(\ell, \hat{y}) \chi\{\mathcal{T}_b^\infty < \mathcal{T}_\blacktriangle^\infty\})$$

since $b < \ell$.

Next, we compare (40) and

(41) $$\sum_{z \in \triangle, \tilde{z}_1 \leq b} \pi(z) h(z) \frac{1}{m(\ell, \hat{y})} \mathrm{E}_z(\mathcal{N}(\ell, \hat{y}) \chi\{\mathcal{T}_b^\infty < \mathcal{T}_\blacktriangle^\infty\}).$$

The difference is less than

$$\sum_{z \in \triangle} \pi(z) h(z) \frac{1}{m(\ell, \hat{y})} \mathrm{E}_z(\mathcal{N}(\ell, \hat{y}) \circ \mathcal{T}_\blacktriangle^\infty \chi\{\mathcal{T}_b^\infty < \mathcal{T}_\blacktriangle^\infty < \infty\})$$

$$\leq \sum_{z \in \Delta, \tilde{z}_1 \leq b} \pi(z) h(z) \sum_{w \in \blacktriangle} P_z(\mathcal{T}_b^\infty < \mathcal{T}_\blacktriangle^\infty < \infty, \mathcal{W}^\infty[\mathcal{T}_\blacktriangle^\infty] = w)$$

$$\times \frac{1}{m(\ell, \hat{y})} \mathrm{E}_w(\mathcal{N}(\ell, \hat{y}))$$

$$\leq C \sum_{z \in \Delta, \tilde{z}_1 \leq b} \pi(z) h(z) \sum_{w \in \blacktriangle} P_z(\mathcal{T}_b^\infty < \mathcal{T}_\blacktriangle^\infty < \infty, \mathcal{W}^\infty[\mathcal{T}_\blacktriangle^\infty] = w)$$

using Lemma 6

$$\leq C \sum_{z \in \Delta, \tilde{z}_1 \leq b} \pi(z) h(z) P_z(\mathcal{T}_b^\infty < \mathcal{T}_\blacktriangle^\infty < \infty)$$

$$\leq C\varepsilon \qquad \text{by (39)}.$$

We therefore investigate the asymptotics of (41). Condition on the point when $\tilde{\mathcal{W}}_1^\infty$ first hits or exceeds $b$. (41) is equal to

$$\sum_{z \in \Delta, \tilde{z}_1 \leq b} \pi(z) h(z) \sum_{0 \leq s} \sum_{v : v_1 = b+s} P_z(\tilde{\mathcal{W}}_1^\infty[\mathcal{T}_b^\infty] = b + s,$$

(42)
$$\hat{\mathcal{W}}^\infty[\mathcal{T}_b^\infty] = \hat{v}, \mathcal{T}_b^\infty < \mathcal{T}_\blacktriangle^\infty)$$

$$\times \frac{1}{m(\ell, \hat{y})} \mathrm{E}_{0,\hat{v}} \mathcal{N}(\ell - (b+s), \hat{y}).$$

However,

$$\frac{1}{m(\ell, \hat{y})} \mathrm{E}_{(0,\hat{v})} \mathcal{N}(\ell - (b+s), \hat{y})$$

$$= \frac{\mathcal{G}((0, \hat{v}), (\ell - (b+s), \hat{y}))}{\mathcal{G}((0, \hat{y}), (\ell, \hat{y}))}$$

$$= \frac{\mathcal{G}((0, \hat{v}), (\ell - (b+s), \hat{y}))}{\mathcal{G}((0, \hat{y}), (\ell - (b+s), \hat{y}))} \frac{\mathcal{G}((0, \hat{y}), (\ell - (b+s), \hat{y}))}{\mathcal{G}((0, \hat{y}), (\ell, \hat{y}))} \to 1$$

by Theorem 1. The first ratio on the right-hand side above is uniformly bounded in $\hat{v}$ as $\ell \to \infty$ by Lemma 6. Moreover,

$$\frac{\mathcal{G}((0, \hat{y}), (\ell - (b+s), \hat{y}))}{\mathcal{G}((0, \hat{y}), (\ell, \hat{y}))} \leq \frac{1}{p((0, \hat{y}), (b+s, \hat{y}))}$$

and this is integrable in $s$ relative to $P_z(\tilde{\mathcal{W}}_1^\infty[\mathcal{T}_b^\infty] = b + s)$, since from assumption A1, the jump sizes and, hence, the excess beyond $b$ are bounded in distribution by an exponential.

Consequently, (42) converges to

$$\sum_{z \in \Delta, \tilde{z}_1 \leq b} \pi(z) h(z) P_z(\mathcal{T}_b^\infty < \mathcal{T}_\blacktriangle^\infty)$$

by dominated convergence. Moreover,

$$\sum_{z \in \triangle, \tilde{z}_1 \leq b} \pi(z) h(z) P_z(\mathcal{T}_b^\infty < \mathcal{T}_\blacktriangle^\infty)$$

(43)
$$\to \sum_{z \in \triangle} \pi(z) h(z) P_z(\mathcal{T}_\blacktriangle^\infty = \infty) = f \quad \text{as } b \to \infty.$$

Equation (38) is within $\varepsilon + C\varepsilon$ of (41) as $\ell$ tends to infinity. Next, (41) tends to $f$ as $\ell$ and $b$ tend to $\infty$. Since $\varepsilon$ is arbitrarily small, we conclude the right-hand side of (37) tends to $\hat{h}^{-1}(\hat{y}) f$, that is,

$$\lim_{\ell \to \infty} e^{\alpha \ell} \frac{1}{m(\ell, \hat{y})} \pi(\ell, \hat{y}) = \hat{h}^{-1}(\hat{y}) f.$$

Finally, we use Theorem 1 to replace $\mathcal{G}((0, \hat{y}), (\ell, \hat{y}))$ by $\varphi(\hat{y})\mathcal{G}((0, \hat{\sigma}), (\ell, \hat{\sigma}))/\varphi(\hat{\sigma}))$ asymptotically. This gives the result. □

6.1. *Asymptotics of two node networks.* Theorem 4 gives the asymptotics of the steady state $\pi$ probability of a discrete time queueing network when the queue at the first node gets big. Below we specialize Theorem 4, using Proposition 1, to prove Theorem 5 for networks with two nodes. We can consider jump processes because, without loss of generality, we can assume the event rate is one, so we can regard the jump rates as transition probabilities of a Markov chain $W$ with kernel $K$ on $S = \{(x, y) | x, y \geq 0, \text{ where } x, y \in \mathbb{Z} \}$ (where $\mathbb{Z}$ denotes the integers). Consequently, $\pi$ is also the stationary distribution of $W$.

Now extend $S$ to $S^\infty = \{(x, y) | y \geq 0, \text{ where } x, y \in \mathbb{Z}\}$ and define $\blacktriangle = \{(x, y) | x \leq 0, y \geq 0\} \cap S^\infty$. The kernel $K$ agrees with a Markov additive kernel $K^\infty$ defined on $S^\infty$ for transitions between points in $S^\infty \setminus \blacktriangle$. $W^\infty$ is the *free* chain with kernel $K^\infty$ on $S^\infty$. The free process $W^\infty$ is a Markov additive process. The additive component $\tilde{W}^\infty$ is the number of customers at the first node, and the Markovian component $\hat{W}^\infty$ is the number of customers at the second node. Markov additive kernels are defined in detail in Section 2, but we can just extend the definition of $K((x, y); (x + u, y + h))$, which is constant for $x > 0$, to $x \in \mathbb{Z}$. We assume $K^\infty$ is nearest neighbor in the Markovian component. We also assume $K^\infty$ is homogeneous, which means $W^\infty$ is a homogeneous random walk off the boundary $\triangle = \{(0, y) : y = 0\}$.

For $y > 0$, define

$$R^+(\theta_1, \theta_2) = \sum_{x', y'} K^\infty((x, y); (x', y')) e^{\theta_1(x'-x)} e^{\theta_2(y'-y)}$$

and for $y = 0$,

$$R^-(\theta_1, \theta_2) = \sum_{x', y'} K^\infty((x, 0); (x', y')) e^{\theta_1(x'-x)} e^{\beta y'}.$$

THEOREM 5. *The bridge case: Let $\theta^b$ be the solution to*

$$\theta_1^b > 0, \tag{44}$$

$$\frac{\partial[R^+(\theta_1^b, \theta_2^b)]}{\partial \theta_2} = 0, \tag{45}$$

$$R^+(\theta^b) = 1. \tag{46}$$

(a) *Suppose $R^-(\theta_1^b, \theta_2^b) \leq R^+(\theta_1^b, \theta_2^b)$.*

(b) *We assume that the distributions $P(V[1] \in \cdot | W^\infty[0] = (0,1), Z[1] = 1)$, $P(V[1] \in \cdot | W^\infty[0] = (0,1), Z[1] = 2)$ and $P(V[1] \in \cdot | W^\infty[0] = (0,1), Z[1] = 0)$ [having z-transforms $S(z)/s$, $P(z)/p$ and $Q(z)/q$ resp.] satisfy the aperiodicity condition given in Proposition 1.*

(c) *We assume cascade paths up the y axis can be neglected by requiring $\sum_{y \geq 0} e^{\theta_2^b y} \pi(0, y) < \infty$ so B4 holds. We also assume B1.*

(d) *We suppose the mean drift of $\mathcal{V}$ above the x-axis is positive; that is, $d_+ > 0$, where $d_+$ is given at (5) and where the constants $\kappa$ and $u = p = q$ are defined as in Section 3, where we define $\mathcal{I}$ as the $h_0$-transform of $K^\infty$ with $h_0(x,y) = \exp(\theta_1^b x) \exp(\theta_2^b y)$. Denote the associated Markov additive process by $\mathcal{W}_0^\infty = (\mathcal{V}_0, \mathcal{Z}_0)$.*

*Then the least action path is the bridge path. If $R^-(\theta_1^b, \theta_2^b) < R^+(\theta_1^b, \theta_2^b)$ (or, equivalently, $\kappa > 0$), then*

$$\pi(\ell, y) \sim f^b \exp(-\theta_1^b \ell) C_+ \ell^{-3/2} \exp(-\theta_2^b y) \varphi(y), \tag{47}$$

*where $f^b$ is a constant obtainable by fast simulation, $\varphi(y) = \frac{p_0}{u}(1 + \kappa y/p_0)^2$ and $C_+$ is given in Proposition 1.*

*If $R^-(\theta_1^b, \theta_2^b) = R^+(\theta_1^b, \theta_2^b)$ (or, equivalently, $\kappa = 0$), then*

$$\pi(\ell, y) \sim f^b C_0 \ell^{-1/2} \exp(-\theta_1^b \ell) \exp(-\theta_2^b y), \tag{48}$$

*where $C_0$ is given in Proposition 1.*

*The jitter case: Suppose we can find $\theta^j = (\theta_1^j, \theta_2^j)$ satisfying $\theta_1^j > 0$ and*

$$R^+(\theta^j) = R^-(\theta^j) = 1. \tag{49}$$

(a) *Suppose $R^-(\theta_1^b, \theta_2^b) > R^+(\theta_1^b, \theta_2^b)$.*

(b) *We assume the aperiodicity of the distribution $P(V[T_0] \in \cdot | W^\infty[0] = (0,0))$ of the additive part $V$ of $W^\infty$ at the time $T_0$ when the Markovian part $Z$ of $W^\infty$ returns to 0. A sufficient condition is given in Lemma 2.*

(c) *We assume $p < q$.*

(d) *We assume cascade paths up the y axis can be neglected by requiring $\sum_{y \geq 0} e^{\theta_2^j y} \pi(0, y) < \infty$.*

(e) *We suppose the mean drift of $\mathcal{V}$ given at C7 in [8] is positive; that is,*

$$\tilde{d}^j = \varphi(0)(S'_0(1) + P'_1(1)) + (1 - \varphi(0))(Q'(1) + S'(1) + P'(1)) > 0,$$

*where $\varphi$ is a probability measure given at C6 in [8],*

(50)
$$\varphi(0) = \Gamma p/p_0 \quad \text{and} \quad \varphi(y) = \Gamma(p/q)^y$$
$$\text{where } \Gamma = (p/p_0 + p/(q-p))^{-1}.$$

*The above constants are defined as follows. Define the harmonic function $h(x,y) = \exp(\theta_1^j x)\exp(\theta_2^j y)$ and construct $\mathcal{W}^\infty$, the h transform of $K^\infty$ having kernel $\mathcal{K}^\infty$. Define $P$, $Q$, $S$, $P_0$ and $S_0$, as well as $p$, $q$, $s$, $p_0$ and $s_0$, using $\mathcal{K}^\infty$ instead of $\mathcal{I}$ as was done near (4).*

*Then we have a jitter case and*

$$\pi(\ell, y) \sim f^j \exp(-\theta_1^j \ell) \frac{1}{\tilde{d}^j} \exp(-\theta_2^j y)\varphi(y)$$

*and $f^j$ is a constant obtainable by fast simulation.*

PROOF. In the bridge case, take $\alpha = \theta_1^b$. We have to check that conditions A0–A7 and B1–B4 in Section 6 hold. These follow almost immediately from our assumptions. The $h_0$-transformed kernel $\mathcal{K}_0^\infty(z,w) := K^\infty(z,w)h_0(w)/h_0(z)$ satisfies the following properties:

(a) $h_0$ is harmonic for the free process off the $x$-axis so $\mathcal{K}_0^\infty((x,y);(x',y'))$ is a probability transition kernel when $y > 0$.

(b) Off the $x$-axis, the mean vertical drift is 0; that is,

$$\sum_w \mathcal{K}_0^\infty((x,y);(x',y'))y' = 0 \quad \text{when } y > 0.$$

(c) We have assumed $\mathcal{K}_0^\infty((x,0),\cdot)$ is substochastic with killing probability $\kappa$.

Let $\hat{\mathcal{K}}_0^\infty$ denote the Markovian part of the transition kernel $\mathcal{K}_0^\infty$. Since the kernel is nearest neighbor, define

(51)
$$p = \hat{\mathcal{K}}_0^\infty(y, y+1), \quad q = \hat{\mathcal{K}}_0^\infty(y, y-1), \quad s = \hat{\mathcal{K}}_0^\infty(y,y)$$
$$\text{for } y > 0,$$
$$p_0 = \hat{\mathcal{K}}_0^\infty(0,1), \quad s_0 = \hat{\mathcal{K}}_0^\infty(0,0), \quad \kappa = 1 - p_0 - s_0.$$

Assumption A5 is automatic since $\hat{\mathcal{K}}^\infty$ has nearest neighbor transitions and it is easy to check that

$$\varphi(0) = 1 \quad \text{and} \quad \varphi(y) = \frac{p_0}{u}a_0^2(y) \quad \text{for } y > 0$$

is the unique (up to multiples) $\sigma$-finite stationary distribution for $\hat{\mathcal{K}}^\infty$. Moreover, the constants are the only harmonic functions for $\hat{\mathcal{K}}^\infty$ since $a_0$ is the unique harmonic function (up to constant multiples) of $\hat{\mathcal{K}}_0^\infty$.

By Section 4, the spectral radius of $\hat{\mathcal{K}}_0^\infty$ is 1. The spectral radius of $\hat{\mathcal{K}}^\infty$ is the same as $\hat{\mathcal{K}}_0^\infty$ and this is 1, so A4 holds. By Proposition 1, if $\kappa > 0$ and $d_+ > 0$, then

$$\mathcal{G}_0((0,0);(\ell,0)) := \sum_{n=0}^\infty (\mathcal{K}_0^\infty)^n((0,0);(0,\ell)) \sim C_+ \ell^{-3/2},$$

where $C_+$ is defined in Proposition 1.

Instead, if $\kappa = 0$ and $d_+ > 0$, then

$$\mathcal{G}_0((0,0);(\ell,0)) \sim C_0 \ell^{-1/2},$$

where $C_0$ is defined in Proposition 1. Moreover, if $\kappa = 0$, then $h_0$ is harmonic for $K^\infty$, so $a_0(y) = 1$ and $\theta^j = \theta^b$.

In the jitter case, take $\alpha = \theta_1^j$. We need to check hypotheses C1–C11 in [8]. These follow immediately from our assumptions. $\varphi$ is as given since the Markovian part is nearest neighbor. $\square$

**7. Modified Jackson networks.** In this section we first give conditions for the stability of a modified Jackson network. We then specialized Theorem 5 to a modified Jackson network with two nodes.

7.1. *Definitions.* Consider a Jackson (1957) network with two nodes. The arrival rate of exogenous customers at Nodes 1 and 2 form Poisson processes with rates $\bar{\lambda}_1$ and $\bar{\lambda}_2$, respectively. The service times are independent, exponentially distributed random variables with mean $1/\mu_1$ and $1/\mu_2$, respectively. Each customer's route through the network forms a Markov chain. A customer completing service at Node 1 is routed to Node 2 with probability $r_{1,2}$ or leaves the system with probability $r_{1,0} := 1 - r_{1,2}$. Routing from Node 2 is defined analogously. So without loss of generality, we are assuming $r_{1,1} = r_{2,2} = 0$. The routing process, service processes and arrival processes are independent.

To ensure that the network is open, we assume that $r_{1,2}r_{2,1} < 1$. Since the network is open, the traffic equations

(52) $$\lambda_i = \bar{\lambda}_i + \lambda_{3-i} r_{3-i,i} \qquad \text{for } i = 1, 2,$$

have a unique solution $(\lambda_1, \lambda_2) = ((\bar{\lambda}_1 + \bar{\lambda}_2 r_{2,1})/(1 - r_{1,2}r_{2,1}), (\bar{\lambda}_2 + \bar{\lambda}_1 r_{1,2})/(1 - r_{1,2}r_{2,1}))$. To eliminate degenerate situations, we assume that $\lambda_1 > 0$ and $\lambda_2 > 0$.

The joint queue length process of this Jackson network forms a Markov process with state space $S = \{0, 1, \ldots\}^2$. Define $\rho_i = \lambda_i/\mu_i$ for $i = 1, 2$. From Jackson (1957), it follows that the stationary distribution for the joint queue length process being in the state $(x,y) \in S$ is $(1-\rho_1)\rho_1^x(1-\rho_2)\rho_2^y$, provided that the stability conditions $\rho_1 < 1$ and $\rho_2 < 1$ hold.

The network that we analyze is a small change from the above network. Suppose that Server 2 has been cross-trained and helps Server 1 whenever Queue 2 is empty. Let $\mu_1^* \geq \mu_1$ be the combined service effort of the two servers at Node 1 when Server 2 is empty. In the section entitled Stability in [9], we proved the modified network will be stable if $\rho_2 < 1$ and $\mu_1^* > (\lambda_1 - \mu_1\rho_2)/(1-\rho_2)$.

We are interested in the rare event of a large deviation in the number of customers at Node 1; that is, more than $\ell$ customers at Node 1 where $\ell$ is large. We will consider choices of the parameters which eliminate the possibility of a large deviation first up the $y$-axis followed by a drift over to $F_\ell$. This means the large deviation paths are along the $x$-axis as was described in Part I. However, if the most likely approach for the Jackson network "jitters" along the $x$-axis, as $\mu_1^*$ increases this approach may eventually become sufficiently difficult so that some other approach becomes most likely. Instead the process travels along the $x$-axis, but instead of jittering along the $x$-axis, the process skims above and only rarely touches the $x$-axis. We refer to this path as a *bridge path*.

The event rate of the modified Jackson network is $\lambda = \bar{\lambda}_1 + \bar{\lambda}_2 + \mu_2 + \mu_1^*$ (since $\mu_1^* > \mu_1$). Without loss of generality, we assume $\lambda = 1$ so we can regard the jump rates as transition probabilities of a Markov chain $W$ with kernel $K$ on $S = \{(x,y) | x, y \geq 0, \text{ where } x, y \in \mathbb{Z}\}$. The modified Jackson network is precisely the homogenization of this chain. $W$ is a nearest neighbor random walk in $S$. Jumps outside of $S$ are suppressed.

7.2. *Asymptotics in the 1-positive recurrent case.* In this section we apply the results in [8] and [14]. Define the harmonic function $h(x,y) = \exp(\theta_1^j x)\exp(\theta_2^j y)$, where $\theta^j$ satisfies (49). We assume $\rho_2 < \exp(-\theta_2^j)$ and $\rho < 1$. The associated Markovian part of the twisted kernel $\hat{\mathcal{K}}^\infty$ is positive recurrent if $\rho < 1$ ($\rho$ is given by Theorem 5). The associated stationary probability $\varphi$ is given in Theorem 5. The only thing to check is

$$\sum_{y=0}^{\infty} \hat{h}(y)\pi(0,y) < \infty.$$

This follows because this sum is bounded by $\sum_{y=0}^{\infty} \exp(\theta_2^j y)\rho_2^y$ by the a priori bound given in Lemma 1 in Part 1. This is finite if $\rho_2 < \exp(-\theta_2^j)$.

In summary:

COROLLARY 1. *If $\rho_2 \leq \exp(-\theta_2^j)$ and if $\rho < 1$ where $\rho = p/q$ and where $p, q$ are defined after* (50), *then the least action path jitters along the $x$-axis*

*and $\varphi$ becomes*

(53) $$\varphi(y) = \begin{cases} \left(1 + \dfrac{r}{1-\rho}\right)^{-1} r\rho^{y-1}, & \text{if } y > 0, \\ \left(1 + \dfrac{r}{1-\rho}\right)^{-1}, & \text{if } y = 0, \end{cases}$$

*where* $r = (\bar{\lambda}_2 e^{\theta_2^j} + \mu_1^* r_{1,2} e^{-\theta_1^j + \theta_2^j})/(\mu_2 r_{2,0} e^{\theta_2^j} + \mu_2 r_{2,1} e^{\theta_1^j - \theta_2^j}).$

7.3. *Asymptotics in the 1-transient case.* In this section we apply the results in Section 6. We assume $\rho_2 < \exp(-\theta_2^b)$, but $\rho > 1$. In this case Proposition 1 gives a solution $\theta^b$ and we take $\alpha$ in A0 to be $\theta_1^b$. A2 holds since the network is nearest neighbor in the additive component, as well as the vertical component. A2 is true by our assumptions on the Jackson network.

$\mu_1^* \geq \mu_1 > 0$ and $\mu_2 > 0$ (or else the network cannot be stable). Moreover, at least one of $\bar{\lambda}_1 > 0$ or $\bar{\lambda}_2 > 0$ or the network is empty. $r_{1,0} > 0$ or $r_{2,0} > 0$ by hypothesis. Finally, $\lambda_1 > 0$ and $\lambda_2 > 0$, which means $\bar{\lambda}_1 + \bar{\lambda}_2 r_{2,1} > 0$ or $\bar{\lambda}_2 + \bar{\lambda}_1 r_{1,2} > 0$. A2 follows immediately.

Next we check hypothesis A2.5. Suppose $\bar{\lambda}_2 > 0$. Next suppose $r_{2,0} > 0$. In the first case when $r_{1,0} > 0$, then from a point $(x,0)$, the chain $W^\infty$ returns to the $x$-axis either at $(x,0)$ by going up and down or at $(x-1,0)$ by going left one. Since we can return to both $x-1$ and to $x$, A2.5 is satisfied. In the other case when $r_{1,2} > 0$, then from a point $(x,0)$ the chain $W^\infty$ returns to the $x$-axis either at $(x,0)$ by going up and down or at $(x-1,0)$ by going northwest one and then down. Again A2.5 is satisfied.

Suppose, on the other hand, that $\bar{\lambda}_2 > 0$ but $r_{2,0} = 0$ so $r_{2,1} > 0$ and $r_{1,0} > 0$. From a point $(x,0)$ the chain $W^\infty$ returns to the $x$-axis either at $(x-1,0)$ taking a step left or to $(x,0)$ by taking a step left, going up one and then going southeast. A2.5 is satisfied. The case when $\bar{\lambda}_2 = 0$ follows in the same way.

Conditions A3–A5 hold automatically and A6 holds since either $\bar{\lambda}_1 > 0$ or $\bar{\lambda}_2 + \mu_2 r_{2,1} > 0$.

Assumptions B2–B3 hold automatically. As for assumption B4, we remark that

$$\sum_{y=0}^{\infty} \hat{h}(y)\pi(0,y) \leq \sum_{y=0}^{\infty} e^{\theta_2^b y} a_0(y) \rho_2^y,$$

so assumption B4 holds if $e^{\theta_2^b}\rho_2 < 1$. Note that this is precisely the condition that appears in Theorem 1 in Part I, ensuring that the least action path is the bridge path or the path jittering along the $x$-axis. Note also that we have appealed to the special structure of the Jackson network, but in a more general case we could use Lyapunov functions.

42R. D. FOLEY AND D. R. MCDONALDThis leaves assumption B1. It suffices to show that the ratio $\mathcal{G}((0,y),(\ell,0))/\mathcal{G}((0,0),(\ell,0))$ is bounded uniformly in $y$ as $\ell \to \infty$ (i.e., we take $\hat{\sigma} = 0$). For paths $\omega'$ of $\mathcal{W}^\infty$ starting from $(0,0)$, define $N_{(\ell,0)}(\omega')$ to be the number of visits by $\mathcal{W}^\infty$ to $(\ell,0)$ following the trajectory $\omega'$. Also define $B$ to be the set of trajectories where the $x$-coordinate $\tilde{\mathcal{W}}_1^\infty$ stays nonnegative, that is, $B = \{\omega' : \tilde{\mathcal{W}}_1^\infty[n] \geq 0, n \geq 0\}$. Similarly, for paths $\omega$ of $\mathcal{W}^\infty$ starting from $(0,y)$, define $N_{(\ell,0)}(\omega)$ to be the number of hits at $(\ell,0)$.4242  R. D. FOLEY AND D. R. MCDONALD42    R. D. FOLEY AND D. R. MCDONALD

This leaves assumption B1. It suffices to show that the ratio $\mathcal{G}((0,y),(\ell,0))/\mathcal{G}((0,0),(\ell,0))$ is bounded uniformly in $y$ as $\ell \to \infty$ (i.e., we take $\hat{\sigma} = 0$). For paths $\omega'$ of $\mathcal{W}^\infty$ starting from $(0,0)$, define $N_{(\ell,0)}(\omega')$ to be the number of visits by $\mathcal{W}^\infty$ to $(\ell,0)$ following the trajectory $\omega'$. Also define $B$ to be the set of trajectories where the $x$-coordinate $\tilde{\mathcal{W}}_1^\infty$ stays nonnegative, that is, $B = \{\omega' : \tilde{\mathcal{W}}_1^\infty[n] \geq 0, n \geq 0\}$. Similarly, for paths $\omega$ of $\mathcal{W}^\infty$ starting from $(0,y)$, define $N_{(\ell,0)}(\omega)$ to be the number of hits at $(\ell,0)$.

Now consider the product space of all paths $\omega$ of $\mathcal{W}^\infty$, which start from $(0,y)$ times paths $\omega'$ which start from $(0,0)$. On this product space we can define a coupled path starting from $(0,0)$, which follows $\omega'$ until $\omega'$ hits the path $\omega$ and then follows $\omega$. Given a path $\omega$ which hits $(\ell,0)$, we note that all paths $\omega' \in B$ must hit the path $\omega$ because the path $\omega'$ is trapped between the $x$-axis and the path $\omega$ and there are no northeast or southwest transitions for Jackson networks.

Define $N_{(\ell,0)}(\omega,\omega')$ to be the number of hits at $(\ell,0)$ by the coupled path. Then, for $\omega' \in B$ and any path $\omega$, $N_{(\ell,0)}(\omega) \leq N_{(\ell,0)}(\omega,\omega')$. Hence,

$$\begin{aligned}
\mathrm{E}_{(0,y)} N_{(\ell,0)}(\omega) \cdot P_{(0,0)}(\omega' \in B) &= \mathrm{E}_{(0,y)} \otimes \mathrm{E}_{(0,0)}(N_{(\ell,0)}(\omega)\chi_B(\omega')) \\
&\leq \mathrm{E}_{(0,y)} \otimes \mathrm{E}_{(0,0)}(N_{(\ell,0)}(\omega,\omega')\chi_B(\omega')) \\
&\leq \mathrm{E}_{(0,y)} \otimes \mathrm{E}_{(0,0)}(N_{(\ell,0)}(\omega,\omega')) \\
&= \mathrm{E}_{(0,0)}(N_{(\ell,0)}(\omega')).
\end{aligned}$$

We conclude $\mathcal{G}((0,y),(\ell,0)/\mathcal{G}((0,0),(\ell,0) \leq P_{(0,0)}(\omega' \in B)^{-1}$; that is, the ratio is uniformly bounded in $\ell$ and $y$:

Remark that the special structure of the modified Jackson network was not really needed to make the above argument work. Bounded jumps are enough.

We can now apply Theorem 5. Since the mean vertical increment is zero, we can simply use the following values:

$$\begin{aligned}
u = p = q &= \bar{\lambda}_2 \exp(\theta_2^b) + \mu_1 r_{1,2} \exp(-\theta_1^b)\exp(\theta_2^b), \\
p_0 &= \bar{\lambda}_2 \exp(\theta_2^b) + \mu_1^* r_{1,2} \exp(-\theta_1^b)\exp(\theta_2^b), \\
s_0 &= \bar{\lambda}_1 \exp(\theta_1^b) + \mu_1^* r_{1,0} \exp(-\theta_1^b)\exp(\theta_2^b), \\
\kappa &= 1 - p_0 - s_0
\end{aligned}$$

and

$$d_+ = -\mu_1 r_{1,2} A_1^{-1} A_2 + \mu_2 r_{2,1} A_1 A_2^{-1} + (\bar{\lambda}_1 A_1 - \mu_1 r_{1,0} A_1^{-1})$$

in the bridge case.

7.4. *Asymptotics in the 1-null recurrent case.* There are practically no changes from the 1-transient case. $\kappa = 0$ if and only if $\theta^b$ also satisfies $R^-(\theta^b) = 0$; that is, if and only if $\theta^b = \theta^j$. In this case $\rho = 1$ and $\varphi(y) = 1$ for all $y$. The application of Theorem 5 is immediate.

**8. A comparison with existing results.** A referee asked us to compare our results with existing results in the literature in [4]. This paper describes two processors with service rates, $\mu_1$ at server one and $\mu_2$ at server two, serving independent streams of Poisson arrivals with rates $\lambda_1$ and $\lambda_2$, respectively. When one server becomes idle it helps the other server and the rate at server one becomes $\mu_1^*$ when server two is idle and the rate at server two becomes $\mu_2^*$ when server one becomes idle. This, nearly a modified Jackson network with $r_{12} = r_{21} = 0$, except that, for simplicity, we did not allow server one to help server two. It is clear that if $\mu_1^*$ is big enough, then the large deviation when the first queue gets big is a bridge.

Unfortunately, we have not been able to derive the asymptotics of the steady state probability of this system from the results in [4]. In the special case where $\mu_1* = \mu_2^* = \mu_1 + \mu_2$ (where $pq = \mu_1\mu_2$ in the notation of [4]), formula (6.4) in [4] should give the asymptotics, but it is not clear how. In the general case when $pq \neq \mu_1\mu_2$, (7.2) in [4] gives a formula for $G(z)$, which is essentially the $z$-transform of $\pi(x,0)$. Again, we could not invert this function.

We did not have better luck with [12], which is a special case of [4] with $\lambda_1 = \lambda_2 = \lambda$, $\mu_1 = \mu_2 = \mu/2$ and $\mu_1^* = \mu_2^* = \mu$. In this, the total number of customers in the system is exactly an $M/M/1$ queue with load $\rho = (\lambda_1 + \lambda_2)/(\mu_1 + \mu_2)$, where the probability the total number of customers is $\ell$ decays like $\rho^\ell$.

Fortunately the analysis in [6] gives the asymptotics of $\pi(\ell,y)$ for the bathroom problem (as it was described in [22]) where couples arrive at a cinema according to a Poisson process with rate $\nu$ and immediately visit the ladies' and men's room. The service rate, at the men's queue is $\alpha$, while the rate at the ladies' queue is $\beta$. This model was extended in [14] by allowing separate arrival streams with rate $\eta$ for single ladies and rate $\lambda$ for single men. We are interested in a large deviation of the men's queue, so let this be the first queue and the ladies' queue the second. In the (unrealistic) case $\alpha < \beta$, the exact asymptotics of $\pi(\ell,y)$ are given in [14]. Omitting the streams of singles (so $\lambda = \eta = 0$),

$$\pi(\ell,y) \sim \frac{f}{\alpha - \nu}\left(\frac{\nu}{\alpha}\right)^\ell \left(1 - \frac{\alpha}{\beta}\right)\left(\frac{\alpha}{\beta}\right)^y,$$

where $f$ is a constant. This is the positive recurrent case and the asymptotics are the same as (7.19) in [6].

The theory in [14] could not handle the (more realistic) case $\alpha > \beta$ (note $\gamma$ in [14] equals $\alpha$ since $\lambda = \eta = 0$). The only result given was that at the stopping time $\tau_\ell$ when the men's queue reaches size $\ell$, the ladies' queue divided by $\ell$ converges to $(\alpha - \beta)/(\alpha - \nu)$. However, the results of this paper do solve this problem because $\beta < \alpha$ gives the 1-transient case. Note that (44), (45), (46) become

$$\theta_1^b > 0, \tag{54}$$

$$\nu e^{\theta_1^b} e^{\theta_2^b} - \beta e^{-\theta_2^b} = 0, \tag{55}$$

$$\nu e^{\theta_1^b} e^{\theta_2^b} + \beta e^{-\beta_2^b} + \alpha e^{-\theta_1^b} = (\nu + \beta + \alpha). \tag{56}$$

Solving (55) gives

$$\exp(\theta_2^b) = \sqrt{\beta/\nu} \exp(-\theta_1^b/2), \tag{57}$$

substituting into (56), and letting $z = \exp(\theta_1^b)$, we get

$$4\beta\nu z^3 = ((\alpha + \beta + \nu)z - \alpha)^2. \tag{58}$$

We recognize this as equation $D_1(z) = 0$ given at (3.2) in [6] (where $\nu = 1$).

To make the connection to [6] (since they only consider the case $\alpha \leq \beta$), we must exchange the labels of the first and second queues so our $\alpha$ becomes $\beta$ and vice versa in [6]. This converts the equation $D_1(z) = 0$ into $D_2(z) = 0$ and according to [6], there is only one solution $z = a_3'$ with $z > 1$, so $a_3' = \exp(\theta_1^b)$. The asymptotics of $\pi(\ell, y)$ are now given by those $p_{y,\ell}$ in (7.20) in [6]. We immediately recognize the asymptotics $\ell^{-3/2}(a_3')^{-\ell}$ for the additive component as was predicted by expression (47).

Next by (57), $\exp(\theta_2^b) = \sqrt{\beta/\nu}(a_3')^{-1/2}$,

$$u = p_0 = p = \nu \exp(\theta_1^b) \exp(\theta_2^b) = \beta \exp(-\theta_2^b), \qquad s = \alpha \exp(-\theta_1^b),$$

$$s_0 = \alpha \exp(-\theta_1^b) + \beta \text{ so } \kappa = 1 - p_0 - s_0 = \beta \exp(-\theta_2^b) - \beta > 0.$$

Hence, the asymptotic distribution of the second queue is given by expression (47):

$$\exp(-\theta_2^b y)\frac{p_0}{u}\left(1 + \frac{\kappa y}{p_0}\right) = (a_3')^{y/2}\left(1 + \left(1 - \sqrt{\frac{\beta}{\nu}}(a_3')^{-1/2}\right)y\right),$$

and this agrees with (7.20) in [6]. Of course, [6] is much more precise since the constant term is given as well.

There is still the null recurrent case to consider. Now $\kappa = 0$ when $\nu e^{\theta_1^b} e^{\theta_2^b} + \alpha e^{-\theta_1^b} + \beta = 1$; that is, when $\beta e^{-\theta_2^b} = \beta$. Hence, $\theta_2^b = 0$, so $\exp(\theta_1^b) = \beta/\nu$, and this implies $\alpha = \beta$ when one considers $\exp(\theta_1^b) = a_3'$ solves $D_2(z) = 0$ in [6]. Hence, Proposition 1 gives the asymptotics $\ell^{-1/2}(\beta/\nu)^{-\ell}$ for the additive component, and this agrees with (7.3) in [6] (since $\alpha = \beta$). Since $\varphi(y) = 1$ in the null recurrent case and since $\theta_2^b = 0$, we conclude

$$\exp(-\theta_2^b y)\frac{p_0}{u}\left(1 + \frac{\kappa y}{p_0}\right) = 1$$

and this agrees with (7.3) in [6].

Of course, we still have to check the conditions A0–A7 and B1–B4, but these are essentially the same as for the modified Jackson network. Even the proof that B1 holds is essentially the same because now we have jumps going northeast, but not southeast or northwest. Condition B4 holds because $\pi(0, y) \leq (1 - \nu/\beta)(\nu/\beta)^y$, so using (57),

$$\sum_{y=0}^{\infty} \pi(0,y)e^{\theta_2^b y} \leq \sum_{y=0}^{\infty}\left(1 - \frac{\nu}{\beta}\right)\left(\frac{\nu}{\beta}\right)^y \left(\sqrt{\frac{\beta}{\nu}}e^{-\theta_1^b/2}\right)^y < \infty.$$

Finally, let us say a word about the matrix-geometric method. The level is the additive component and the phase is the Markovian component of a Markov additive process with kernel $K^\infty$. For a QBD process, the additive component is nearest neighbor. Denote the transition probabilities from phase $i$ at level $\ell$ to phase $j$ at level $\ell - 1$, respectively, $\ell$ and $\ell + 1$ by $A$, respectively, $B$ and $C$. Hence, we can represent the Feynman–Kac kernel as

$$\hat{K}_\gamma(i,j) = \sum_t e^{\gamma t}K((i,0);(j,t)) = (\nu^{-1}A + B + \nu C)_{ij}$$

if $\nu = \exp(\gamma)$. Solving the Riccatti equation associated with the matrix-geometric method is equivalent to finding a Perron–Frobenius eigenvector of the Feynman–Kac kernel with eigenvalue one. That means solving $(\nu^{-1}A + B + \nu C)\mathbf{y} = \mathbf{y}$ is equivalent to finding $\gamma$ and $\hat{h} = \mathbf{y}$ so $\hat{K}_\gamma \hat{h} = \hat{h}$. Equivalently, the solution to the Riccatti equation gives a harmonic function $h(x,i) = \exp(\gamma x)\hat{h}(i)$ for $K^\infty$.

**Acknowledgments.** David R. McDonald wishes to thank the members of the Department of Mathematics of the Indian Institute of Science in Bangalore for their hospitality while working on this paper. He also thanks François Baccelli and the École Normale Supérieure for their support while the revision was completed. Finally, we thank Bina Bhattacharyya for her help.

## REFERENCES

[1] ADAN, I. J. B. F., WESSELS, J. and ZIJM, W. H. M. (1993). A compensation approach for 2-D Markov processes. *Adv. in Appl. Probab.* **25** 783–817. MR1241929

[2] CHUNG, K. L. (1967). *Markov Chains with Stationary Transition Probabilities*, 2nd ed. Springer, New York. MR217872

[3] DAVIS, B. and MCDONALD, D. (1995). An elementary proof of the local central limit theorem. *J. Theoret. Probab.* **8** 693–701. MR1340834

[4] FAYOLLE, G. and IASNOGORODSKI, R. (1979). Two coupled processors: The reduction to a Riemann-Hilbert problem. *Z. Wahrsch. Verw. Gebiete* **47** 325–351. MR525314

[5] FELLER, W. (1971). *An Introduction to Probability Theory and Its Applications* **II**. Wiley, New York. MR270403

[6] FLATTO, L. and HAHN, S. (1984). Two parallel queues created by arrivals with two demands I. *SIAM J. Appl. Math.* **44** 1041–1053. MR759714

[7] FLATTO, L. and MCKEAN, H. P. (1977). Two queues in parallel. *Comm. Pure Appl. Math.* **30** 255–263. MR464426

[8] FOLEY, R. and MCDONALD, D. (2001). Join the shortest queue: Stability and exact asymptotics. *Ann. Appl. Probab.* **11** 569–607. MR1865017

[9] FOLEY, R. and MCDONALD, D. (2005). Large deviations of a modified Jackson network: Stability and rough asymptotics. *Ann. Appl. Probab.* **15** 519–541. MR2114981

[10] IGNATIOUK-ROBERT, I. (2001). Sample path large deviations and convergence parameters. *Ann. Appl. Probab.* **11** 1292–1329. MR1878299

[11] KESTEN, H. (1995). A ratio limit theorem for (sub) Markov chains on $\{1, 2, \ldots\}$ with bounded jumps. *Ann. Appl. Probab.* **27** 652–691. MR1341881

[12] KONHEIM, A., MEILIJSON, I. and MELKMAN, A. (1981). Processor-sharing of two parallel lines. *J. Appl. Probab.* **18** 952–956. MR633243

[13] MCDONALD, D. (1979). On local limit theorems for integer valued random variables. *Theory Probab. Appl.* **24** 613–619. MR541375

[14] MCDONALD, D. (1999). Asymptotics of first passage times for random walk in a quadrant. *Ann. Appl. Probab.* **9** 110–145. MR1682592

[15] MEYN, S. P. and TWEEDIE, R. L. (1993). *Markov Chains and Stochastic Stability*. Springer, New York. MR1287609

[16] NEUTS, M. F. (1981). *Matrix-Geometric Solutions in Stochastic Models*; *an Algorithmic Approach*. Johns Hopkins Univ. Press. MR618123

[17] NEY, P. and NUMMELIN, E. (1987a). Markov additive processes I. Eigenvalue properties and limit theorems. *Ann. Probab.* **15** 561–592. MR885131

[18] NEY, P. and SPITZER, F. (1966). The Martin boundary for random walk. *Trans. Amer. Math. Soc.* **121** 116–132. MR195151

[19] OLVER, F. W. (1974). *Asymptotics and Special Functions*. Academic Press, New York. MR435697

[20] SENETA, E. and VERE-JONES, D. (1966). On quasi-stationary distributions in discrete-time Markov chains with a denumerable infinity of states. *J. Appl. Probab.* **3** 403–434. MR207047

[21] SERFOZO, R. (1999) *Introduction to Stochastic Networks*. Springer, New York. MR1704237

[22] SHWARTZ, A. and WEISS, A. (1994). *Large Deviations for Performance Analysis*. Chapman and Hall, New York.

[23] TAKAHASHI, Y., FUJIMOTO, K. and MAKIMOTO, N. (2001). Geometric decay of the steady-state probabilities in a quasi-birth-and-death process with a countable number of phases. *Stoch. Models* **17** 1–24. MR1852862

[24] WOESS, W. (2000). *Random Walks on Infinite Graphs and Groups*. Cambridge Univ. Press. MR1743100

School of Industrial
  and Systems Engineering
Georgia Institute of Technology
765 Ferst Drive
Atlanta, Georgia 30332-0205
USA

Department of Mathematics
  and Statistics
University of Ottawa
Ottawa, Ontario
Canada K1N 6N5
e-mail: dmdsg@uottawa.ca